\documentclass[12pt]{amsart}     

\usepackage{tkz-fct}

\usepackage[margin=1.15in]{geometry}

\usepackage[pagebackref=true,
]{hyperref} 
\hypersetup{colorlinks=true,linkcolor=blue,filecolor=magenta,citecolor=green,urlcolor=blue,final} 

 \usepackage{amssymb}

\usepackage{amsthm,amssymb,amsmath,amscd}
\usepackage{amsfonts,color}
\usepackage{latexsym}
\usepackage[all]{xy}
\usepackage{ulem}

\usepackage{mathptmx}

\newtheorem{thm}{Theorem}[section]
\newtheorem{lem}[thm]{Lemma}
\newtheorem{cor}[thm]{Corollary}
\newtheorem{prop}[thm]{Proposition}
\theoremstyle{definition}
\newtheorem{df}[thm]{Definition}
\theoremstyle{remark}
\newtheorem{rem}[thm]{Remark}

\newtheorem{example}[thm]{Example}

\numberwithin{equation}{section}

\newtheorem*{ack}{Acknowledgements}

\newcommand{\bC}{{\mathbb C}}
\newcommand{\bD}{{\mathbb D}}

\newcommand{\bL}{{\mathbb L}}

\newcommand{\bP}{{\mathbb P}}
\newcommand{\bR}{{\mathbb R}}
\newcommand{\bQ}{{\mathbb Q}}
\newcommand{\bZ}{{\mathbb Z}}
\newcommand{\bT}{{\mathbb T}}

\newcommand{\cD}{{\mathcal D}}

\newcommand{\cF}{{\mathcal F}}

\newcommand{\cH}{{\mathcal H}}

\newcommand{\cL}{{\mathcal L}}
\newcommand{\cM}{{\mathcal M}}

\newcommand{\cO}{{\mathcal O}}

\newcommand{\cS}{{\mathcal S}}

\newcommand{\wti}{\widetilde}

\newcommand{\Gr}{\text{\rm Gr}}
\newcommand{\DR}{\hbox{\rm DR}}

\newcommand{\lra}{\longrightarrow}

\newcommand{\Int}{{\rm Int}}
\newcommand{\Relint}{{\rm Relint}}

\newcommand{\MHM}{{\rm MHM}}
\newcommand{\MHS}{{\rm MHS}}
\newcommand{\rat}{{\it {rat}}}

\def\sig{\sigma}
\def\Sig{\Sigma}

\def\be{\begin{equation}}
\def\ee{\end{equation}}
\def\bt{\begin{thm}}
\def\et{\end{thm}}
\def\bc{\begin{cor}}
\def\ec{\end{cor}}
\def\br{\begin{rem}}
\def\er{\end{rem}}
\def\bp{\begin{prop}}
\def\ep{\end{prop}}
\def\bl{\begin{lem}}
\def\el{\end{lem}}
\def\bn{\begin{enumerate}}
\def\en{\end{enumerate}}
\def\bex{\begin{example}}
\def\eex{\end{example}}
\def\bd{\begin{df}}
\def\ed{\end{df}}

\begin{document}                        


\title[]{Weighted Ehrhart theory  via mixed Hodge modules \\ on toric varieties}

\author[L. Maxim ]{Lauren\c{t}iu Maxim}
\address{L. Maxim : Department of Mathematics, University of Wisconsin-Madison, 480 Lincoln Drive, Madison WI 53706-1388, USA, \newline
{\text and} \newline Institute of Mathematics of the Romanian Academy, P.O. Box 1-764, 70700 Bucharest, ROMANIA.}
\email {maxim@math.wisc.edu}

\author[J. Sch\"urmann ]{J\"org Sch\"urmann}
\address{J.  Sch\"urmann : Mathematische Institut,
          Universit\"at M\"unster,
          Einsteinstr. 62, 48149 M\"unster,
          Germany.}
\email {jschuerm@math.uni-muenster.de}

\subjclass[2020]{14M25, 14C17, 14C40, 52B20}
\keywords{Toric varieties, lattice polytopes, lattice points, Ehrhart polynomial, reciprocity, Hirzebruch classes, mixed Hodge module, intersection cohomology}

\date{\today}

\begin{abstract}
We give a cohomological and geometrical interpretation for the weighted Ehrhart theory of a full-dimensional lattice polytope $P$, with Laurent polynomial weights of geometric origin. For this purpose, 
we calculate the motivic Chern and Hirzebruch characteristic classes of a mixed Hodge module complex $\cM$ whose underlying cohomology sheaves are constant on the $\bT$-orbits of the toric variety $X_P$ associated to $P$. Besides motivic coefficients, this also applies to the intersection cohomology Hodge module. 
We introduce a corresponding generalized Hodge $\chi_y$-polynomial of the ample divisor $D_P$  on $X_P$. Motivic properties of these characteristic classes are used to express this Hodge polynomial in terms of a very general weighed lattice point counting and the corresponding weighted Ehrhart theory. We introduce, for such a mixed Hodge modules complex $\cM$ on $X$, an Ehrhart polynomial $E_{P,\cM}$ generalizing the Hodge polynomial of $\cM$ and satisfying a reciprocity formula and a purity formula fitting with the duality for mixed Hodge modules. This Ehrhart polynomial and  its properties depend only on a Laurent polynomial weight function on the faces $Q$ of $P$. In the special case of the intersection cohomology mixed Hodge module, the weight function corresponds to Stanley's $g$-function  of the polar polytope of $P$, hence it depends only on the combinatorics of $P$. In particular, we obtain a combinatorial formula for the intersection cohomology signature.
\end{abstract}

\maketitle   

\section{Introduction}
We give a cohomological and geometrical perspective for the Ehrhart theory of a full-dimensional lattice polytope $P$, expressed here in the form of a weighted Ehrhart theory, with Laurent polynomial weights $f_Q(y)\in \bZ[y^{\pm 1}]$ of geometric origin being attached to the faces $\emptyset \neq Q \preceq P$ of the polytope (with the partial order given by face inclusion). The classical combinatorial approach to reciprocity can be linearly extended to this weighted Ehrhart theory. Motivation for using such weights comes from our prior work \cite{MS} on Hirzebruch characteristic classes of toric varieties (associated to such polytopes), as well as from the use of Stanley's $g$-functions in the recent generalized Ehrhart reciprocity theorem of \cite{BGM}.

\subsection{Combinatorial and geometrical posets and their identification} We review here some basic notions explaining the relation between lattice polytopes and toric geometry; for more details see \cite{CLS, Fu2}.

Let $M\cong \bZ^n$ be a lattice and $P \subset M_{\bR}:=M\otimes \bR \cong \bR^n$ be a full-dimensional lattice polytope. 
To $P$ one associates a fan $\Sig=\Sig_P \subset N_\bR:=N\otimes \bR$, with $N$ the dual lattice of $M$, called the {inner normal fan} of $P$.  To the fan $\Sig_P$ one further associates a projective toric variety $X=X_P$ with torus $\bT$ (whose character lattice is $M$) and an ample Cartier divisor $D=D_P$. There is a order-reversing one-to-one correspondence between the faces $\emptyset \neq Q \preceq P$ of $P$ and the cones $\sig_Q$ of $\Sig_P$, with $\dim_\bR(Q)=n-\dim_\bR(\sig_Q)$. There is an order-reversing one-to-one correspondence between the cones $\sig_Q$ of $\Sig_P$ and the $\bT$-orbits $O_Q$ of the toric variety $X_P$, with $\dim_\bC(O_Q)=n-\dim_\bR(\sig_Q)$. Then the identification between the faces $\emptyset \neq Q \preceq P$ of $P$ and the torus orbits $O_Q$ of the toric variety $X_P$ is order-preserving, with $\dim_\bR(Q)=\dim_\bC(O_Q)$.

If, in addition, $P$ contains the origin in its interior, then its {\it polar polytope} $P^\circ\subset N_\bR$
is defined as in \cite[Section 1.5]{Fu2}, as a full-dimensional lattice polytope with respect to the lattice $N$, and containing the origin in the interior. By taking cones at the origin of $N_\bR$ over the proper faces of $P^\circ$, with $\emptyset$ corresponding to the origin, one gets the same lattice fan $\Sig_P$ (hence the same toric variety $X_P$). This correspondence between the proper faces $\emptyset \preceq Q^\circ \prec P^\circ$ of $P^\circ$ and the cones of $\Sig_P$ is order-preserving, increasing real dimension by $1$. (Here $\dim_\bR \emptyset =-1$ by convention.) This induces an order-reversing one-to-one correspondence between the faces $Q$ of $P$, and  the faces $Q^\circ$ of the polar polytope $P^\circ$, switching the roles of polytopes and emptysets seen as faces. Moreover,  for a proper face $\emptyset \neq Q \prec P$, one has $\dim_\bR(Q)+\dim_\bR(Q^\circ)=n-1$.


The original polytope $P \subset M_\bR$ is classically used for counting lattice points (and the corresponding Ehrhart theory), fitting with sections of the ample line bundle $\cO_{X_P}(D_P)$ on $X_P$ \cite{Dan,Fu2,CLS}. The torus orbits of $X_P$ give a natural $\bT$-invariant Whitney stratification, which is particularly useful in our geometric interpretation of weighted lattice points counting via the theory of Hirzebruch homology classes $T_{y*}$ of \cite{BSY} and the corresponding Hirzebruch-Riemann-Roch theorem; see \cite{MS}. The polar polytope $P^\circ$ appears in Stanley's work \cite{St} for recursively defining his {\it $g$-polynomials} $g_{Q^\circ}(t) \in \bZ[t]$, with $g(\emptyset)=1$, for the Eulerian graded poset given by the faces $\emptyset \preceq Q^\circ \preceq P^\circ$. It is well-known \cite{St, BM,DL,F,Sa} that these $g$-polynomials are related to the intersection cohomology complex of $X_P$. In the use of the $g$-polynomials, we will (implicitly) assume that $P$ contains the origin in its interior (so that the polar polytope is well-defined).

\subsection{Weighted Ehrhart theory} Let us now assign to the above posets Laurent polynomial weights via a weight vector  $f=\{f_Q\}$, with $f_Q(y) \in \bZ[y^{\pm 1}]$ indexed here by the faces $\emptyset \prec Q \preceq P$ of $P$. We assemble these weights, together with the combinatorics of $P$, into the {\it weighted Ehrhart polynomial} $E_{P,f}(\ell, y)$, as follows.
\bd For $\ell \in \bZ_{>0}$, define the weighted Ehrhart polynomial of $P$ and $f$ by
\be\label{ep1i} E_{P,f}(\ell, y):=\sum_{Q \preceq P} f_Q(y) \cdot (1+y)^{\dim(Q)} \cdot | \Relint({{\ell}} Q) \cap M | ,\ee
with $|-|$ denoting the cardinality of a finite set, and 
$ \Relint({{\ell}} Q)$ is the relative interior of the face ${\ell} Q$ of the dilated polytope $\ell P$.
\ed

By classical Ehrhart theory (e.g., see \cite{Dan, CLS}), $E_{P,f}(\ell, y)$ has the following properties:
\begin{enumerate}
\item $E_{P,f}(\ell, y)$ is obtained by evaluating a polynomial $E_{P,f}(z, y)$ at $z=\ell \in \bZ_{>0}$.
\item (Constant term) \ For $\ell=0$, $$E_{P,f}(0, y)=\sum_{Q \preceq P} f_Q(y) \cdot (-1-y)^{\dim(Q)},$$ i.e., evaluating $| \Relint({{\ell}} Q) \cap M |$ at $\ell=0$ as $(-1)^{\dim(Q)}$.
\item (Reciprocity formula) \ For $\ell \in \bZ_{>0}$,
\be\label{rei} E_{P,f}(-\ell, y)=\sum_{Q \preceq P} f_Q(y) \cdot (-1-y)^{\dim(Q)} \cdot | \ell Q \cap M |.\ee
\end{enumerate}

Assume now that the origin is an interior point of $P$, so that we can consider the weight vector given by Stanley's $g$-polynomials \be f_Q(y)=g_{Q^\circ}(-y)=:\wti{g}_Q(-y)\ee for the polar polytope of $P$, with $g_{\emptyset}(-y)=\wti{g}_P(-y)=1$. In this setup, it was shown in \cite[Theorems 1.3 and 2.6]{BGM} that the following {\it purity} property holds:
\be\label{pi}E_{P,f}(-\ell, y)=(-y)^n \cdot E_{P,f}(\ell, 1/y).\ee
This is in fact a special case of the the results in loc.cit., for the constant polynomial $\varphi=1$. Moreover, our  definition \eqref{ep1i} for $E_{P,f}(\ell, y)$ equals $G_{\varphi=1}(\ell,y)$ from \cite[formula (14)]{BGM}. In the case when $P$ is a simple polytope, the polar polytope $P^\circ$ is simplicial, so that $g_{Q^\circ}(-y)=1$, for all faces $Q$ of $P$. In this case, as explained in \cite{BGM}, the equality \eqref{pi} implies the Dehn-Sommerville relations for $P$.

\subsection{Geometric viewpoint on weighted Ehrhart theory}
The condition that $P$ is simple translates into the fact that the corresponding toric variety $X_P$ is a rational homology manifold, i.e., the constant sheaf on $X_P$ agrees with the shifted intersection cohomology complex: $$Q_{X_P}\simeq IC_{X_P}[-n]=:IC'_{X_P}.$$
Then the classical Dehn-Sommerville relations correspond to the Poincar\'e duality isomorphism for rational cohomology, motivating the use of intersection cohomology coefficients for non-simple lattice polytopes. 

In this paper, we give a geometric interpretation and proofs of the above properties.
For instance, using the intersection Hodge module ${IC'}^H_{X_P}$ and its purity, we give a geometric proof of formula \eqref{pi}, as it will be explained below.
For generalizations to the equivariant context, which will imply these properties for weighted Ehrhart polynomials depending on any homogeneous polynomial functions $\varphi:M_\bR \to \bC$ as in \cite{BGM}, see the authors' paper \cite{MS3}.

\subsubsection{Geometric realization of weight vectors and duality}
First, we identify the weight vectors $\{f_Q\}$, seen as $\bZ[y^{\pm 1}]$-valued functions on the set of faces $Q$ of $P$, with the $\bT$-invariant $\bZ[y^{\pm 1}]$-valued constructible functions on $X=X_P$, via $$\{f_Q\} \mapsto \sum_{Q \preceq P} 1_{O_Q} \cdot f_Q \in F^{\bT}(X)[y^{\pm 1}].$$
Note that $F^{\bT}(X)[y^{\pm 1}]$ is a $\bZ[y^{\pm 1}]$-module freely generated by the $\bT$-invariant constructible functions $1_{O_Q}\in F^{\bT}(X) \subset F^{\bT}(X)[y^{\pm 1}].$

There is a tautological injective homomorphism
$$\tau: F^{\bT}(X)[y^{\pm 1}] \lra K_0(var/X)[\bL^{-1}]$$
to the Grothendieck group of varieties over $X$, localized at the Lefschetz motive $\bL=[\bC \to pt]$, given by $$1_{O_Q} \cdot f_Q(y)  \mapsto [O_Q \hookrightarrow X] \cdot  f_Q(-\bL).$$ The image $\tau(F^{\bT}(X)[y^{\pm 1}])$ of this monomorphism is the $\bZ[\bL^{\pm 1}]$-submodule freely generated by the images of the classes $[O_Q \hookrightarrow X]\in K_0(var/X) \to K_0(var/X)[\bL^{-1}]$.

Most of our calculations can already be done on this motivic level. But in order to identify the intersection cohomology complex and duality we need to go one step further, to the Grothendieck group $K_0(\MHM(X))$ of Saito's mixed Hodge modules on $X$. There is a canonical homomorphism (see \eqref{chih})
$$\chi_{Hdg}: K_0(var/X)[\bL^{-1}] \to K_0(\MHM(X))$$
commuting with duality (introduced on the left hand side in \cite{Bit}, and with the duality $\bD_X$ on the right hand side given by the Verdier duality of mixed Hodge modules), see \cite{Sc} for details. We have
$$\chi_{Hdg}([O_Q \hookrightarrow X] \cdot  f_Q(-\bL))=[(j_Q)_!\bQ^H_{O_Q}] \cdot f_Q(-[\bQ(-1)]),$$
with $j_Q:O_Q \hookrightarrow X$ the orbit inclusion and $\bQ^H_{O_Q}$ denoting the constant Hodge module on $O_Q$. The restriction of $\chi_{Hdg}$ to the subgroup $\tau(F^{\bT}(X)[y^{\pm 1}])$ is injective, with image given by the free $\bZ[u^{\pm 1}]$-submodule of $K_0(\MHM(X))$ generated by elements of the form $[(j_Q)_!\bQ^H_{O_Q}]$, where $u=[\bQ(-1)]\in K_0(\MHM(pt))\simeq K_0(\MHS^p)$. Here the module structure on $K_0(\MHM(X))$ is induced from the exterior product for the identification $X\simeq pt \times X$, and $\MHS^p$ is the abelian category of graded polarizable $\bQ$-mixed Hodge structures, with $\bQ(-1)$ the corresponding Tate Hodge structure of weight $2$.

The ambient Grothendieck groups $K_0(var/X)[\bL^{-1}]$ and $K_0(\MHM(X))$ are only used formally in the background. For our applications to the toric context, we only need to work with the (sub)groups 
\be\label{gh} \{f_Q\} \longleftrightarrow F^{\bT}(X)[y^{\pm 1}] \longleftrightarrow  \tau(F^{\bT}(X)[y^{\pm 1}]) 
\longleftrightarrow  \chi_{Hdg}\left(\tau(F^{\bT}(X)[y^{\pm 1}])\right)\ee
realizing our weight vectors in the language of constructible functions, motivic, and mixed Hodge module contexts, respectively. These identifications are linear over the following ring coefficient isomorphisms, with $-y=\bL=u$:
$$\bZ[y^{\pm 1}]  \longleftrightarrow \bZ[y^{\pm 1}] \longleftrightarrow \bZ[\bL^{\pm 1}] \longleftrightarrow \bZ[u^{\pm 1}],$$
identifying the coordinate basis elements as follows: $$ e_Q:=\{\delta_{Q,Q'}\}  \longleftrightarrow 1_{O_Q}  \longleftrightarrow [O_Q \hookrightarrow X]  \longleftrightarrow [(j_Q)_!\bQ^H_{O_Q}].$$

By their definitions, all groups homomorphisms in \eqref{gh} are surjective. In our geometric applications, we work in the group $\chi_{Hdg}\left(\tau(F^{\bT}(X)\right)$, recovering the corresponding weight vector $\{f_Q\}$ via $$f_Q(y):=\chi_y(i_{x_Q}^*[\cM])$$ as the Hodge $\chi_y$-polynomial of the stalk class $i_{x_Q}^*[\cM]$, with $[\cM]\in \chi_{Hdg}\left(\tau(F^{\bT}(X)\right)$, and $i_{x_Q}:\{x_Q\} \hookrightarrow X$ the inclusion of a point $x_Q$ chosen in the orbit $O_Q$. Then our aim is to show that many invariants of $[\cM]$ depend only of this choice of a weight function.

Specific to the toric context is that the last submodule $\chi_{Hdg}\left(\tau(F^{\bT}(X)[y^{\pm 1}])\right)$ is preserved under duality, since each torus orbit has a $\bT$-invariant affine neighborhood of product type. Since ${IC'}^H_{X}$ is a pure Hodge module of weight $n$, so that $\bD_X{IC'}^H_{X}\simeq {IC'}^H_{X}[2n](n)$, we get the following purity property:
\be\label{dua} \bD_X [{IC'}^H_{X}]=u^{-n} \cdot 
[{IC'}^H_{X}].
\ee
Note that this duality has the following module property: $\bD_X(u \cdot -)=u^{-1}\cdot \bD_X(-)$.
This induces a duality involution $\bD$ on the $\bZ[y^{\pm1}]$-module of weight functions with the corresponding 
module property $\bD(y \cdot -)=y^{-1} \cdot \bD(-)$.
Under this identification, $[{IC'}^H_{X}]$ corresponds in our language to the collection of Stanley's $g$-functions $\wti{g}:=\{\wti{g}_Q(-y)\}$, i.e., $\wti{g}_Q(-y)=\chi_y(i_{x_Q}^*[{IC'}^H_{X}])$, see \cite{BM,F,DL,Sa}. 

\subsubsection{Degree calculations} The power of working with motivic and mixed Hodge module Grothendieck groups resides in their functorial calculus, e.g., pushforward for a proper map commuting with cross-product (e.g., module structure), and duality. Applying this to the constant map $X \to pt$, together with the Hodge $\chi_y$-polynomial homomorphism, we get for example:
$$ [O_Q \hookrightarrow X] \mapsto \chi_y([O_Q \hookrightarrow pt])=(-y-1)^{\dim(Q)}.$$
When translated back to the side of weight vectors, this degree map corresponds to the constant term $E_{P,f}(0, y)$ of the weighted Ehrhart polynomial, also motivating the presence of the factor $(y+1)^{\dim(Q)}$ in the definition of $E_{P,f}(\ell, y)$. Note that this degree map on the constructible function side $F^{\bT}(X)[y^{\pm 1}]$ is not the usual one, which one recovers just for $y=-1$ (with the convention $0^0=1$).

The advantage of working on the mixed Hodge module side is the use of mixed Hodge module complexes $\cM\in D^b\MHM(X)$ and their important properties to represent weight vectors. Pushing $[\cM]$ down to a point yields $[H^*(X;\cM)]\in K_0(\MHS^p)$, with cohomology groups getting induced (graded polarizable) mixed Hodge structures. Here the corresponding Hodge $\chi_y$-polynomial is defined in terms of the Hodge filtration. 

As an example, $\cM={IC'}^H_X$ is pure of weight zero in the sense of Saito \cite{Sa1,Sa2}, so that 
$IH^i(X):=H^i(X;{IC'}^H_X)$ is a pure Hodge structure of weight $i$. On the other hand, $[IH^*(X)]\in \bZ[u^{\pm 1}]$, i.e., it is of Tate type. These two properties imply the vanishing of odd intersection cohomology groups, so that $I\chi_y(X):=\chi_y([IH^*(X)])$ coincides with the intersection cohomology Poincar\'e polynomial $IP_X(t)=h_P(t^2)$, up to the change of variable $t^2=-y$. As is well-known, one also has a similar purity and odd vanishing property for the stalks of $[{IC'}^H_X]$ (e.g., see \cite{Sa}). So we get 
\be\begin{split}
h_P(t^2)&=E_{P,f}(0, -t^2)\\
&=\sum_{Q \preceq P} \wti{g}_Q(t^2) \cdot (t^2-1)^{\dim(Q)}\\
&=\sum_{\emptyset \preceq Q^\circ \prec P^\circ} {g}_{Q^\circ}(t^2) \cdot (t^2-1)^{n-1-\dim(Q^\circ)},
\end{split}
\end{equation}
which gives a geometric interpretation of Stanley's recursion for the $h$-polynomial (see also \cite{F,DL,Sa}).
Furthermore, using the Hodge index theorem for intersection cohomology \cite{MSS}, one gets for $y=1$ the following combinatorial formula for the Goresky-MacPherson signature of $X=X_P$:
\be\label{sigci} sign(X)=E_{P,\wti{g}}(0, 1)=\sum_{Q \preceq P} \wti{g}_Q(-1) \cdot (-2)^{\dim(Q)}.\ee
Lastly, pushing \eqref{dua} to a point and applying the $\chi_y$-homomorphism yields the duality formula $$I\chi_y(X)=(-y)^n \cdot I\chi_{1/y}(X).$$ This already recovers Stanley's ``master duality'' theorem \cite[Thm. 3.14.9]{St2}, as well as the generalized Dehn-Sommerville relations for the $h$-polynomial of a lattice polytope (boiling down to Poincar\'e duality for $IH^*(X)$). In terms of the corresponding Ehrhart polynomials, this duality formula reads as
$$E_{P,\wti{g}}(0, y)=(-y)^n \cdot E_{P,\wti{g}}(0, 1/y),$$
which proves \eqref{pi} for $\ell=0$.

\subsubsection{Generalized Hirzebruch-Riemann-Roch and weighted Ehrhart theory}
We can now explain our geometric interpretation of the weighted Ehrhart theory, in terms of the Hirzebruch class transformation $$T_{y*}:K_0(\MHM(X)) \to H_{2*}(X;\bQ)[y^{\pm 1}]$$ of \cite{BSY,Sc}. Our approach extends, and reduces to, Danilov's classical work \cite{Dan}, by using two key formulae for the Hirzebruch classes of the basic elements $[(j_Q)_!\bQ^H_{O_Q}]$ and their duals $[(j_Q)_*\bQ^H_{O_Q}]$; see  \eqref{drt} and \eqref{drtc}. This Hirzebruch class transformation commutes with proper pushforwards, cross-products, and with duality in the following sense:
\be\label{dui} T_{y*}\circ \bD_X=\bD \circ T_{y*},\ee
where $\bD_X$ is induced from the duality of mixed Hodge modules, and $\bD$ on $H_{2i}(X;\bQ)$ is defined by $(-1)^i \cdot id$,  extended to $H_{2*}(X;\bQ)[y^{\pm}]$ by $y \mapsto 1/y$. On a point space, $T_{y*}$ reduces to the Hodge $\chi_y$-polynomial. The Hirzebruch classes $T_{y*}([\cM])$ fit into the following {\it generalized Hirzebruch-Riemann-Roch formula} (see Theorem \ref{gHRR})
\be
\chi_y(X,D;\cM)=\int_{X} ch(\cO_X(D)) \cap T_{y*}([\cM]),
\ee
with $D=D_P$ the ample Cartier divisor on $X=X_P$ associated to $P$; see Section \ref{sec4} for the precise definition of $\chi_y(X,D;\cM)$.

Then for $[\cM]\in   \chi_{Hdg}\left(\tau(F^{\bT}(X)[y^{\pm 1}])\right)$, and $\ell \in \bZ_{>0}$, we prove in Proposition \ref{wlpca} the following identity:
\be\label{wcai}
\begin{split}
\chi_y(X,\ell D;\cM) &=
\sum_{k=0}^{n} \left(\frac{1}{k!} \int_{X} [D]^k \cap T_{y*}([\cM])  \right) \ell^k \\
 &=\sum_{Q \preceq P} \chi_y(i_{x_Q}^*[\cM]) \cdot (1+y)^{\dim(Q)} \cdot | \Relint({{\ell}} Q) \cap M | \\
 &=E_{P,f}(\ell,y), 
\end{split}
\ee
with weight vector $f_Q(y)=\chi_y(i_{x_Q}^*[\cM])$, 
where the summation on the right is over the faces $Q$ of $P$, and $x_Q\in O_{\sigma_Q}$ is a point in the orbit of $X_P$ associated to (the cone $\sigma_Q$ in the inner normal fan of $P$ corresponding to) the face $Q$. This reproves in cohomological terms the {\it polynomial property} for $E_{P,f}(\ell,y)$, together with the formula for the constant term $$E_{P,f}(0,y)=\chi_y([H^*(X;\cM)]).$$

Using the above cohomological description of $E_{P,f}(\ell,y)$, together with the duality property \eqref{dui}, we get the following identity (see Theorem \ref{repr}):
\be\label{r1i}
\begin{split}
E_{P,\cM}(-\ell, y)=E_{P, \bD_X\cM}(\ell, {1}/{y}),
\end{split}
\ee
which yields the {\it reciprocity formula} \eqref{rei}, 
for the weight vector defined as above from the stalk classes of $\cM$.
If $\cM$ is such a self-dual pure Hodge module of weight $n$ on $X=X_P$, so that $\bD_X [\cM]=u^{-n} \cdot [\cM]$ with $u=[\bQ(-1)]$,  
then the following {\it purity property} follows:
\be\label{r1bi}
E_{P,\cM}(-\ell, y)=(-y)^n \cdot E_{P, \cM}(\ell, {1}/{y}).
\ee
In particular, if $\cM={IC'}_X^H$, corresponding to the weight function $f_Q(y)=\wti{g}_Q(-y)$, the purity property \eqref{r1bi}, together with classical Ehrhart reciprocity, {\it proves} formula \eqref{pi}. 

\subsection{Structure of the paper} The paper is structured as follows. 

In Section \ref{sec2}, we give a general introduction to the $K$-theoretic Hodge-Chern class $\DR_y$, which is used for the definition of the homology Hirzebruch classes $T_{y*}$ of mixed Hodge module complexes. We explain here the basic calculus and properties of these classes. 

Preparing the ground for the toric context, in Section \ref{sec3} we detail the stratum-wise calculation of these classes.  Under the technical assumption of (virtual) stratum-wise constancy of the underlying cohomology sheaves of the coefficients, this reduces to a motivic type calculation. These results apply to more coefficients that those from the subgroup $\chi_{Hdg}\left(\tau(F^{\bT}(X)[y^{\pm 1}])\right)$.

In Section \ref{sec4}, for a Cartier divisor $D$ on a compact complex algebraic variety $X$, we define for a mixed Hodge module complex $\cM$ a generalized Hodge polynomial $\chi_y(X,D;\cM)$, which can be computed via a generalized Hirzebruch-Riemann-Roch theorem in terms of the Chern character of $\cO_X(D)$ and the Hirzebruch class $T_{y*}([\cM])$ of the Grothendieck class of $\cM$. For $D=0$ this reduces to the usual Hodge polynomial $\chi_y(X;\cM)$ of $\cM$. We also prove in this section a key duality formula  \eqref{duf} for $\chi_y(X,D;\cM)$.

In Section \ref{sec5} we specialize to the toric context, and make the characteristic class formulae from the previous sections much more explicit. For a first read, those interested in applications to Ehrhart theory can start with this section, taking for granted a few technical results from the prior sections.

In the final Section \ref{sec6}, we expand on the weighted Ehrhart theory, giving a geometric view and proofs of the results discussed in this introduction.

\smallskip

\begin{ack}
L. Maxim 
acknowledges support from the project ``Singularities and Applications'' - CF 132/31.07.2023 funded by the European Union - NextGenerationEU - through Romania's National Recovery and Resilience Plan.
J. Sch\"urmann was funded by the Deutsche Forschungsgemeinschaft (DFG, German Research Foundation) Project-ID 427320536 -- SFB 1442, as well as under Germany's Excellence Strategy EXC 2044 390685587, Mathematics M\"unster: Dynamics -- Geometry -- Structure. L. Maxim and  J. Sch\"urmann also thank the Isaac Newton Institute for Mathematical Sciences for the support and hospitality during the program ``Equivariant methods in geometry'' when work on this paper was undertaken.
\end{ack}

\section{Homology Hirzebruch classes via mixed Hodge modules}\label{sec2}
Let $X$ be a complex algebraic variety, and denote by 
$\MHM(X)$ the abelian category of algebraic mixed Hodge modules on $X$, with $K_0(\MHM(X)) \simeq  K_0(D^b\MHM(X))$ the Grothendieck group of (bounded complexes of) mixed Hodge modules on $X$. Recall that $\MHM(pt)\simeq \MHS^p$ is (equivalent to) the category of (graded-)polarizable $\bQ$-mixed Hodge structures.

For any such variety $X$, Saito constructed a functor of triangulated categories 
$$Gr^F_p\DR: D^b\MHM(X) \lra D^b_{\rm coh}(X) $$
commuting with proper pushforwards,
with $Gr^F_p\DR(\cM)\simeq 0$ for almost all $p$ and $\cM$ fixed; see  \cite[Sec.2.3]{Sa1},  \cite[Sec.3.10, Prop.3.11]{Sa2} and \cite[Sec.1]{Sa3}. Here $D^b_{\rm coh}(X)$ denotes the bounded derived category of $\cO_X$-modules with coherent cohomology.

\bex If $X$ is a closed subvariety of a  $n$-dimensional complex algebraic manifold $Z$, and $\cM\in \MHM_X(Z)$ is a mixed Hodge module on $Z$ with support on $X$, then
$Gr^F_p\DR(\cM)$ is the complex associated to the de Rham complex $\DR(\cM)$ of the underlying algebraic left $\cD_Z$-module $\cM$ with its integrable connection $\nabla$:
$$ \DR(\cM)=[\cM \overset{ \nabla}{\lra}  \cM\otimes_{\cO_Z} \Omega^1_Z  \overset{ \nabla}{\lra} \cdots \overset{ \nabla}{\lra} \cM\otimes_{\cO_Z} \Omega^n_Z]
$$
with $\cM$ in degree $-n$, filtered  by
$$
F_p \DR(\cM) =[F_p\cM \overset{ \nabla}{\lra} F_{p+1}\cM \overset{ \nabla}{\lra} \cdots  \overset{ \nabla}{\lra} F_{p+n}\cM\otimes_{\cO_Z} \Omega^n_Z]. 
$$
By properties of mixed Hodge modules, the graded parts of the De Rham complex belong to $D^b_{coh}(X)$.
In the (ambient) smooth case, we use the same symbol for both the mixed Hodge module and the underlying filtered (left) $\cD_Z$-module
 \qed
\eex

For any complex algebraic variety $X$, the transformations $Gr^F_p\DR$ induce group homomorphisms 
$$Gr^F_p\DR: K_0(\MHM(X)) \lra K_0(D^b_{\rm coh}(X)) \simeq K_0(X),$$
with $K_0(X)$ the Grothendieck group of the abelian category of coherent $\cO_X$-modules.
\bd[Brasselet--Sch\"urmann--Yokura \cite{BSY,Sc}]
The {\it Hodge--Chern class transformation} of a complex algebraic variety $X$ is:
$$ \DR_y:K_0(\MHM(X)) \lra K_0(X) \otimes \bZ[y^{\pm 1}] $$
\begin{equation*}
\begin{split} 
\DR_y([\cM])&:=\sum_{i,p}\,(-1)^i \big[\cH^i Gr^F_{-p}\DR(\cM)\big] \cdot (-y)^p \\
&=\sum_{p}\, \big[Gr^F_{-p}\DR(\cM)\big] \cdot (-y)^p
\end{split}
\end{equation*}
The {\it Hirzebruch class transformation} is defined by
$$ T_{y*}:=td_* \circ \DR_y :K_0(\MHM(X)) \to H_*(X)[y^{\pm 1}] $$
with $td_*:K_0(X) \to H_*(X)$ the Todd class transformation \cite{BFM} of the singular Grothendieck--Riemann--Roch theorem of Baum-Fulton-MacPherson, linearly extended over $\bZ[y^{\pm 1}]$, and $H_*(X):=H_{2*}^{BM}(X)\otimes \bQ$ is the even degree Borel-Moore homology of $X$ with $\bQ$-coefficients.
\ed

\br
As shown in \cite{BSY,Sc}, the transformations $\DR_y$ and (by Riemann-Roch) $T_{y*}$ 
commute with proper pushforward and with cross-products $\boxtimes$. 
\er

\bex[Degree]
Let us illustrate the above definition when $X=pt$ is a point space. In this case we have
$$K_0(pt)=\bZ,\quad H_*(pt)=\bQ,\quad \MHM(pt)=\MHS^p,$$
where $\MHS^p$ is, as before,  Deligne's category of graded-polarizable mixed $\bQ$-Hodge
structures, with switching the increasing $\cD$-module filtration to a decreasing Hodge filtration so that on a point space this identification gives $Gr^F_{-p}\DR=Gr^p_F$.
By definition, we then have for $H^{\bullet}\in D^b\MHS^p$ that
\be\label{316} 
\chi_y([H^{\bullet}]):=
DR_{y}([H^{\bullet}])=T_{y*}([H^{\bullet}])=
\sum_{j,p}\,(-1)^j\dim_{\bC}\Gr_F^pH_{\bC}^j\,(-y)^p\ee
is the corresponding Hodge $\chi_y$-polynomial, defining a ring homomorphism 
$$\chi_y:K_0(\MHS^p) \lra \bZ[y^{\pm 1}]\subset \bQ[y^{\pm 1}].$$
E.g., for $n \in \bZ$, 
$\chi_y(\bQ(n))=(-y)^{-n}$, where $\bQ(n)$ denotes the Tate Hodge structure of weight $-2n$ on the vector space $\bQ$. 
Similarly, $$\chi_y(X):=\chi_y([H^\bullet_c(X;\bQ_X^H)])$$ is the (compactly supported) Hodge polynomial of $X$. 
E.g., $\chi_y((\bC^*)^n)=(-y-1)^n$.
\qed
\eex

\bex By functoriality, 
for $X$ compact, $\cM\in D^b\MHM(X)$, and $k:X \to pt$ the constant map to a point, 
\be\label{eq2} \chi_y(X;\cM):=\chi_y([H^\bullet(X;\cM)])=\chi_y([k_*\cM])=\int_X T_{y*}([\cM])
\ee
is the degree of the Hirzebruch class of $\cM$. For $\cM=\bQ_X^H$ the constant Hodge module, this fits with $\chi_y(X)$ as defined above. \qed
\eex

\bd
Let $K_0(\MHS^p)^{\rm Tate}$ be the subring of $K_0(\MHS^p)$ generated by $u:=[\bQ(-1)]$ and $u^{-1}:=[\bQ(1)]$, with a surjective ring homomorphism $\bZ[u^{\pm 1}] \to K_0(\MHS^p)^{\rm Tate}$. This is an isomorphism since its composition with $\chi_y$ is the isomorphism $\bZ[u^{\pm 1}] \to \bZ[y^{\pm 1}], \ u \mapsto -y$.
\ed

Let $K_0(var/X)$ be the motivic {relative Grothendieck group} of complex algebraic varieties over $X$,
i.e., the free abelian group generated by isomorphism classes $[f]:=[f: Y\to X]$ of morphisms $f$ to $X$,
divided out be the usual {scissor relation}.  The  pushdown $f_!$, cross-product $\boxtimes$ and pullback $g^*$
for these relative Grothendieck groups 
are defined by composition, cross-product and, resp., pullback of arrows. These transformations are linear over $K_0(var/pt)$ by the cross-product for $pt \times X\simeq X$.
Denote by $\bL:=[\bC \to pt]$ the Lefschetz motive.
There is a natural group homomorphism on the localized Grothendieck group
\be\label{chih} \chi_{Hdg}: K_0(var/X)[\bL^{-1}] \to K_0(\MHM(X)), \ \ \ [f:Y\to X]\mapsto [f_!\bQ^H_Y]\:,\ee
mapping $\bL$ to $u=[\bQ(-1)] \in K_0(\MHS^p)=K_0(\MHM(pt)$. This transformation commutes with pushdown $f_!$, cross-product $\boxtimes$ and pullback $g^*$. 
The compositions
$$mC_y:=DR_y \circ \chi_{Hdg}: K_0(var/X)[\bL^{-1}] \to K_0(X)[y]$$
and
$$T_{y*}:=T_{y*} \circ \chi_{Hdg}: K_0(var/X)[\bL^{-1}] \to H_*(X)[y]$$
are the motivic Chern and Hirzebruch classes from \cite{BSY,Sc}, as studied in the toric context in \cite{MS} (without localization).

As usual, for a variety $X$ we set 
$$mC_y(X):=mC_y([id_X])=\DR_y([\bQ_X^H])$$
 and
$$T_{y*}(X):=T_{y*}([id_X])=T_{y*}([\bQ_X^H])$$
for $\bQ_X^H$ the constant Hodge module (complex) on $X$.
Note that these classes can also be defined in terms of the filtered Du Bois complex (see \cite{BSY}). 

Moreover, if $X$ is pure-dimensional we set
$$IC_y(X):=\DR_{y}([{IC'}_X^H]), \ \ \ 
IT_{y*}(X):=T_{y*}([{IC'}_X^H]),$$
with  ${IC'}_X^H:=IC_X^H[-\dim(X)]$. If $X$ is compact, the degree of $IT_{y*}(X)$ is the $I\chi_y$-polynomial $$I\chi_y(X):=\chi_y(X;{IC'}_X^H).$$ If $X$ is moreover a rational homology manifold (e.g., a complex algebraic $V$-manifold like a simplicial toric variety), then ${IC'}_X^H \simeq \bQ_X^H$ (which also motivates our shift convention).

\br[Normalization] If $X$ is {smooth} and $\cM=\bQ_X^H$ is the constant Hodge module, then 
$$Gr^F_{-p}\DR(\bQ_X^H)\simeq \Omega^p_X[-p]$$
for $0 \leq p \leq \dim(X)$, and it is $0$ otherwise.  
So, $$mC_y(X)=\sum_{p=0}^{\dim(X)} [\Omega^p_X]\cdot y^p=:\Lambda_y(\Omega^1_X),$$ the total $\Lambda$-class of the cotangent bundle of $X$. If, moreover, $X$ is compact, the degree of $mC_y(X)$ is exactly the Hirzebruch $\chi_y$-genus of $X$.
\er

\bex
Assume $X_\Sig$ is a toric variety (or a complex algebraic $V$-manifold). Then
\be\label{zero}
mC_y(X_\Sig)=\sum_{p = 0}^{\dim(X)} [\widehat{\Omega}_{X_\Sig}^p] \cdot y^p \in K_0(X_\Sig)[y], 
\ee
where  $\widehat{\Omega}_{X_\Sig}^p$ denotes the corresponding sheaf of Zariski differential $p$-forms. In fact, this formula holds even for a torus-invariant closed algebraic subset $X:=X_{\Sig'}\subseteq X_{\Sig}$ (i.e., a closed union of torus-orbits) corresponding to a star-closed subset $\Sig' \subseteq \Sig$, with the sheaves of differential $p$-forms  $\wti{\Omega}^p_{X}$ as introduced by Ishida \cite{I}:
\be\label{zeroi}
mC_y(X)=\sum_{p = 0}^{\dim(X)} [\widetilde{\Omega}_{X}^p] \cdot y^p \in K_0(X)[y].
\ee
This follows, as in \cite{MS}, from \cite[Proposition 4.2]{I}, via the identification of these Ishida sheaves with the (shifted) graded parts of the filtered Du Bois complex (and similarly for the Du Bois complex of a complex algebraic $V$-manifold, by \cite[Section 2.7]{MSS}). \qed
\eex

\br\label{fic}
For later use, note also that ${IC'}_X^H$ is a direct summand of $f_*\bQ^H_M \in D^b\MHM(X)$ for any resolution of singularities $f:M\to X$. By functoriality, this yields that $Gr^F_{-p}\DR({IC'}_X^H)$ is a summand of $Rf_*\Omega^p_M[-p] \in D^b_{\rm coh}(X)$, therefore $Gr^F_{-p}\DR({IC'}_X^H) \simeq 0$ for all $p<0$ or $p>\dim(X)$. In particular, $I\chi_y(X) \in \bZ[y]$ is a polynomial in $y$, with $I\chi_0(X)$ an intersection cohomology version of the arithmetic genus of $X$.
\er

\br[Effect of duality]\label{r9}
As shown in \cite[Corollary 5.19, Remark 5.20]{Sc}, the transformations $\DR_y$ and $T_{y*}$  commute with duality, i.e., $$\DR_y\circ \bD_X=\bD \circ \DR_y, \ \ \ T_{y*}\circ \bD_X=\bD \circ T_{y*},$$
where $\bD_X$ is induced from the duality of mixed Hodge modules, and $\bD$ on $K_0(X)$ is induced by Grothendieck duality (e.g., see \cite[Part 1, Sect.7]{FM}, \cite[Part 1]{LH}) extended to $K_0(X)[y^{\pm}]$ by $y \mapsto 1/y$. Similarly, one uses the duality involution $^\vee$  in homology, given on $H_i(X)$ by $(-1)^i \cdot id$,  extended to  $H_*(X)[y^{\pm}]$ by $y \mapsto 1/y$. 
The relation between the two duality formulas above is via a corresponding duality formula for the Todd class transformation:
\be
td_* \circ \bD = \bD \circ td_*
\ee
Note also that the effect of the Tate twist $(n)$ on $\DR_y$ and $T_{y*}$ is just multiplication by $(-y)^{-n}$, i.e.,  for $\cM \in D^b\MHM(X)$ with $\cM(n):=\cM \boxtimes \bQ(n)$, one has
$$\DR_y([\cM(n)])=(-y)^{-n} \cdot \DR_y([\cM]),$$
and similarly for $T_{y*}$. A similar statement holds for the motivic transformations $mC_y$ and $T_{y*}$, with duality on the localization $K_0(var/X)[\bL^{-1}]$ defined 
as in \cite[Section 4B]{Sc} and \cite{Bit}. All these duality transformations are functorial for proper maps, so that for $X$ compact with $k:X \to pt$ proper, one gets
\be\label{eq5}\chi_y(X;\bD_X \cM)=\chi_{1/y}(X;\cM).\ee
Indeed, over a point space, Grothendieck duality and homological duality are trivial, while the mixed Hodge module duality is the usual duality of mixed Hodge structures, e.g., $\bD_{pt}\bQ(n)=\bQ(-n)$.
\er

\bex 

Applying the above remark to the self-dual object $IC_X^H$, with $\bD_X IC_X^H \simeq IC_X^H(n)$, for a pure-dimensional algebraic variety $X$ with $n=\dim(X)$, one gets the identity \be\label{dual} IT_{y*}(X)=(-y)^n \cdot \left(\sum_{i=0}^n (-1)^i \cdot IT_{1/y,i}(X)  \right) .\ee
In particular, if $y=1$, then 
\be\label{dual1} IT_{1*}(X)=(-1)^n \cdot \left(\sum_{i=0}^n (-1)^i \cdot IT_{1,i}(X) \right)=(-1)^n \cdot \left(IT_{1*}(X) \right)^\vee  ,\ee
i.e., $IT_{1,i}(X)=0$ if $n+i$ is odd. 
For $X$ compact, by taking degrees and using Remark \ref{fic}, one gets
\be I\chi_y(X)=(-y)^n \cdot  I\chi_{1/y}(X) \in \bZ[y]\subset \bZ[y^{\pm 1}],\ee
so, if $y=1$ and $n$ odd, then $I\chi_1(X)=0$.\qed
\eex

\br
In the case of a projective complex algebraic variety $X$, one has by the intersection cohomology Hodge index theorem \cite[Sect.3.6]{MSS2} the equality
\be\label{hin} I\chi_1(X)=sign(X),\ee
with $sign(X)$ the Goresky-MacPherson intersection homology signature. 
If $X$ is moreover a rational homology manifold  (e.g., a complex algebraic $V$-manifold like a simplicial toric variety), $sign(X)$ coincides with the usual signature defined via the intersection pairing in rational homology.
\er


\section{Computations via stratifications}\label{sec3}

Recall the following result, which in the case of toric varieties will be applied to the Whitney stratification given by the torus orbits:

\bp\cite[Prop.5.1.2]{MSS} \label{p7}
For a complex variety $X$, fix $\cM \in D^b\MHM(X)$ with underlying bounded constructible complex $K:=\rat(\cM)\in D^b_c(X)$.
Let $\cS=\{S\}$ be a complex algebraic stratification of $X$ so
that for any stratum $S\in\cS$, $S$ is smooth, ${\bar S}\setminus S$ is a union of
strata, and the sheaves $L_{S,\ell}:=\cH^{\ell}K|_S$ are local systems on $S$ for any $\ell \in \bZ$.
If $$j_S:S\overset{i_{S,{\bar S}}}{\hookrightarrow} {\bar S} \overset{i_{{\bar S},X}}{\hookrightarrow} X$$ is the inclusion map of a stratum $S\in \cS$, then: 
\be\label{f1} [\cM] 
=\sum_{S,\ell}\,(-1)^{\ell}\,\big[(j_S)_!L_{S,\ell}^H \big]
=\sum_{S,\ell}\,(-1)^{\ell}\,(i_{{\bar S},X})_*\big[(i_{S,{\bar S}})_!L_{S,\ell}^H \big] \in K_0(\MHM(X)), 
\ee
where $L_{S,\ell}^H=H^{\ell + \dim (S)} (j_S)^* \cM[-\dim(S)]$ is the shifted smooth mixed Hodge module on $S$ with $L_{S,\ell}=rat(L_{S,\ell}^H)$. 
In particular,
\be\label{f2}
 \DR_y([\cM])=\sum_{S,\ell}\,(-1)^{\ell}\,(i_{{\bar S},X})_*\DR_y\big([(i_{S,{\bar S}})_!L_{S,\ell}^H ]\big), \ee
and \be\label{f3} 
 T_{y*}([\cM])  =\sum_{S,\ell}\,(-1)^{\ell}\,(i_{{\bar S},X})_*T_{y*}\big([(i_{S,{\bar S}})_!L_{S,\ell}^H ]\big).\ee
\ep

\bex[Constant cohomology sheaves along strata]\label{ccs}
If the above stratification $\cS$ of $X$ can be chosen so that the local systems $L_{S,\ell}$ on $S$ are actually {\it constant} for any $S$ and $\ell$, with stalk (the mixed Hodge structure) $L_{S, s, \ell}$ for $s\in S$ fixed, then it follows by {\it rigidity} (e.g., see \cite[Section 3.1]{CMS}) that the (admissible) variation of mixed Hodge structures $L^H_{S,\ell}$ on $S$ (corresponding to a smooth mixed Hodge module as in the above Proposition) is the constant variation. This implies that, if $a_S:S \to pt$ is the constant map and $i_{s,S}:\{s\} \hookrightarrow S$ denotes the inclusion, then 
$$L^H_{S,\ell}\simeq a_S^* i_{s,S}^* L^H_{S,\ell} = a_S^* L^H_{S, s, \ell}\simeq \bQ_S^H \otimes a_S^* L^H_{S, s, \ell} \simeq \bQ_S^H \boxtimes L^H_{S, s, \ell}.$$
Since for any variety $Z$, $K_0(\MHM(Z)$ is a unitary $K_0(\MHM(pt)$-module (via pullback and tensor product),  under our assumptions formula \eqref{f1} becomes:
\be\label{f4}
\begin{split}
 [\cM]
&=\sum_{S,\ell}\,(-1)^{\ell}\,(i_{{\bar S},X})_*\big[(i_{S,{\bar S}})_!\bQ_S^H \big] \cdot [L^H_{S, s, \ell}]\\
&=\sum_{S}\,(i_{{\bar S},X})_*\big[(i_{S,{\bar S}})_!\bQ_S^H \big] \cdot [\cM_s],
\end{split}
\ee
with $\cM_s:=i_s^*\cM$ and $i_s:=j_S \circ i_{s,S} : \{s\} \hookrightarrow X$ for $s\in S$ chosen, so that
\be\label{f4b}
[\cM]=\sum_{S}\,(i_{{\bar S},X})_*\big[(i_{S,{\bar S}})_!\bQ_S^H \big] \cdot i_s^*[\cM].
\ee
In particular, upon applying $\DR_y$ and using the fact that this transformation commutes with cross-products, we get
\be\label{f5}
\begin{split}
 \DR_y([\cM]) 
& =\sum_{S}\,(i_{{\bar S},X})_*\DR_y\big[(i_{S,{\bar S}})_!\bQ_S^H \big] \cdot \chi_y(i_s^*[\cM]) \\
&=\sum_{S}\,(i_{{\bar S},X})_*mC_y ([S \hookrightarrow {\bar S}])  \cdot \chi_y(i_s^*[\cM]),
 \end{split}
 \ee
 with $\chi_y([\cM_s])\in \bZ[y^{\pm 1}]$ regarded as a weight (depending on $\cM$) attached to the stratum containing $s$.
 A similar formula can be obtained for $T_{y*}([\cM])$ by further applying $td_*$ to formula \eqref{f5}.\qed
\eex

\bd\label{sc}
We call $\cM \in D^b\MHM(X)$ a {\it stratum-wise constant} complex of mixed Hodge modules (with respect to an algebraic stratification $\cS$ as above)  if the cohomology sheaves of the underlying constructible complex $K=rat(\cM)\in D^b_c(X)$ are constant along the strata of the stratification. More generally, we call a class $[\cM]\in K_0(\MHM(X))$ {\it virtually (Tate-type) stratum-wise constant} if equality \eqref{f4b} holds (and $i_s^*[\cM] \in K_0(\MHS^p)^{\rm Tate}$) .
\ed

\br
Note that by the above formulae \eqref{f4} and \eqref{f4b}, a stratum-wise constant complex of mixed Hodge modules $\cM$ represents a virtually stratum-wise constant class $[\cM]$. Only the virtually stratum-wise constant property was needed to derive formula \eqref{f5}. For example, if, in the notations of Example \ref{ccs}, all the local systems $L_{S,\ell}$ on $S$ underlie {\it unipotent} variations of mixed Hodge structures (i.e., the graded parts $Gr_k^W L_{S,\ell}$ of the weight filtration are constant, for all $k \in \bZ$), then $[\cM]\in K_0(\MHM(X))$ is virtually stratum-wise constant.
\er

\bex\label{exvcc}
In the notations of Example \ref{ccs}, let $[L^H] \in K_0(\MHM(pt))$ be the Grothendieck class of a mixed Hodge structure $L^H$. Then, for a fixed stratum $S \in \cS$ of $X$, the class
$$[L^H] \cdot (i_{{\bar S},X})_*\big[(i_{S,{\bar S}})_!\bQ_S^H \big]=\big[(i_{{\bar S},X})_*(i_{S,{\bar S}})_! L_S^H]=[\cM_S],$$ where $L^H_S=\bQ_S^H  \boxtimes L^H$, is virtually stratum-wise constant, with $i^*_s[\cM_S]=[L^H]$ for $s\in S$ all other stalks equal to zero. This shows that a class $[\cM]$ is virtually stratum-wise (Tate-type) constant if, and only if, it belongs to the (free) $K_0(\MHM(pt))$-submodule (resp., $K_0(\MHS^p)^{\rm Tate}$-submodule) of $K_0(\MHM(X))$ generated by classes $(i_{{\bar S},X})_*\big[(i_{S,{\bar S}})_!\bQ_S^H \big]$, $S \in \cS$. 
\qed
\eex

\br
Under the assumption of (virtual) stratum-wise constancy, formula \eqref{f5} reduces the calculation of Hodge-Chern and Hirzebruch classes to the motivic calculus of \cite{MS}, up to assigning the above mentioned Laurent polynomial weights $ \chi_y(i_s^*[\cM])$ to each stratum $S$, with $s \in S$ chosen. For instance, the stratum-wise constancy assumption is automatically satisfied in the toric context for $\bT$-equivariant mixed Hodge modules, e.g., the $IC$-complex, see \cite[Lemma 1.2]{Ta}.
\er

\bex[Mapping situation]
Let $f:Y \to X$ be a proper morphism of algebraic varieties.
Let $\cS=\{S\}$ be a complex algebraic stratification of $X$ so
that for any stratum $S\in\cS$, $S$ is smooth, and ${\bar S}\setminus S$ is a union of
strata. Assume that all sheaves $R^\ell f_*\bQ_Y|_S$ are constant on $S$ for any $\ell \in \bZ$ and $S \in \cS$. Then in the above notations, one gets by proper base change,
\be\label{fmap}
\begin{split}
f_*mC_y(Y)&=mC_y([f])=\DR_y([f_*\bQ_Y^H]) \\
&=\sum_{S}\,(i_{{\bar S},X})_*mC_y ([S \hookrightarrow {\bar S}])  \cdot \chi_y(f^{-1}(s)),
\end{split}
\ee
for $s \in S$ chosen. If $X$ (and hence $Y$) is compact, by taking degrees in \eqref{fmap}, we get
$$\chi_y(Y)=\sum_{S} \chi_y({S})  \cdot \chi_y(f^{-1}(s))
=\sum_{S} \left( \chi_y({\bar S}) - \chi_y(\partial {\bar S}) \right)  \cdot \chi_y(f^{-1}(s)),$$ with $\partial {\bar S}= {\bar S}\setminus S$. For similar {\it stratified multiplicative formulae}, see \cite{CLS}.

For instance, if $f$ is a closed inclusion of a union of strata, then only strata in $Y$ contribute in formula \eqref{fmap}, each with multiplicity $1$. Another example is provided by a projective bundle $\bP(V) \to X$, with $V \to X$ an algebraic vector bundle of constant rank $r+1$, so that all multiplicities are equal, and given by $\chi_y(\bP^r)$. \qed
\eex

\br
More generally, the property of stratum-wise constant higher direct image sheaves $R^\ell f_*\bQ_Y|_S$ is satisfied if $f:f^{-1}(S) \to S$ is a Zariski locally trivial fibration for all strata $S \in \cS$. This property also holds for a $\bT$-equivariant proper morphism $f:Y \to X$, with $X$ a toric variety with torus $\bT$ (by \cite[Lemma 1.2]{Ta}).
\er

The terms $mC_y ([S \hookrightarrow {\bar S}]) =\DR_y\big[(i_{S,{\bar S}})_!\bQ_{S}^H \big]$ appearing in formula \eqref{f5} are computed via resolutions of singularities as follows (see \cite[Prop.5.2.1]{MSS} and \cite[Prop.2.2]{MS}):

\bp\label{p8}
Let $i_{S,Z}:S\hookrightarrow Z$ be a smooth partial compactification of a stratum $S$ so that 
$D:=Z\setminus S$ is a simple normal crossing divisor, and $i_{S,{\bar S}}=\pi_Z\circ i_{S,Z}$ for a proper morphism
$\pi_Z:Z\to{\bar S}$. Then we have
\be
mC_y ([S \hookrightarrow {\bar S}]) = \DR_y(\big[ (i_{S,{\bar S}})_!\bQ^H_S\big])=
(\pi_Z)_*\big[ \mathcal{O}_Z(-D)  \otimes \Lambda_y \Omega_{Z}^1(\log(D)\big] \in K_0({\bar S})[y],
\ee
and
\be
\begin{split} T_{y*}\big([(i_{S,{\bar S}})_!\bQ^H_S]\big) &=
\sum_{q\geq 0}(\pi_Z)_*td_{*}\big[\cO_Z(-D) \otimes\Omega_Z^q(\log D)\big]y^{q}\\
&=
(\pi_Z)_*td_{*}\big[\cO_Z(-D) \otimes \Lambda_y\Omega_Z^q(\log D)\big] \in H_*({\bar S})[y].
\end{split}
\ee
\ep

\br\label{ps}
For later use, let us also mention here the following formula, using the above notations (see, e.g., \cite[Example 5.8]{Sc}):
\be
\DR_y(\big[ (i_{S,{\bar S}})_*\bQ^H_S\big])=
(\pi_Z)_*\big[ \Lambda_y \Omega_{Z}^1(\log(D)\big] \in K_0({\bar S})[y],
\ee
and similarly for $T_{y*}$.
\er

These formulae become more explicit and much simplified in the toric context, as it will be discussed in Section \ref{sec5} below.


\section{A generalized HRR for arbitrary coefficients}\label{sec4}
Let $X$ be a compact complex algebraic variety of dimension $n$, with $D$
a Cartier divisor  on $X$. Fix $\cM\in D^b\MHM(X)$. Note that, for any integer $p$, we have
$$Gr^F_{p}\DR(\cM) \otimes \cO_X(D) \in D^b_{\rm coh}(X),$$ 
and we define the {\it Hodge polynomial of $(X,D,\cM)$} by the formula
\be\label{new}
\begin{split}
\chi_y(X,D;\cM)&:=\sum _{p\in \bZ} \chi\left(X,Gr^F_{-p}\DR(\cM) \otimes \cO_X(D) \right) \cdot (-y)^p\\
&=: \chi(X, \DR_y([\cM]) \otimes \cO_X(D)),
\end{split}
\ee
which only depends on $[\cM] \in K_0(\MHM(X))$.
In the last notation, we extend the tensor product linearly over $\bZ[y^{\pm 1}]$.
This is well defined since $Gr^F_p\DR(\cM)\simeq 0$ for almost all $p$. 

\bex If $X$ is smooth and $\cM=\bQ_X^H$ is the constant Hodge module, then 
$$Gr^F_{-p}\DR(\bQ_X^H)\simeq \Omega^p_X[-p]$$
for $0 \leq p \leq \dim(X)$, and it is $0$ otherwise, 
so $$\chi_y(X,D;\bQ_X^H)=\sum _{p=0}^n \chi\left(X, \Omega^p_X \otimes \cO_X(D) \right) \cdot y^p,$$
wich is exactly $\chi_y(X,\cO_X(D))$, the Hirzebruch polynomial of $D$, see \cite[Sect.15.5]{H}.\qed
\eex

\bex
Assume $X$ is a toric variety or a complex algebraic $V$-manifold. Then
$$\chi_y(X,D;\bQ_X^H)=\sum _{p=0}^n \chi\left(X, \widehat{\Omega}^p_X \otimes \cO_X(D) \right) \cdot y^p,$$
where  $\widehat{\Omega}_{X}^p$ denotes the corresponding sheaf of Zariski differential $p$-forms. As before, a similar formula holds for a torus-invariant closed algebraic subset $X_{\Sig'}\subseteq X_{\Sig}$ of a toric variety $X_\Sig$ (i.e., a closed union of torus-orbits) corresponding to a star-closed subset $\Sig' \subseteq \Sig$, with the sheaves of differential $p$-forms  $\wti{\Omega}^p_{X_{\Sig'}}$ as introduced by Ishida \cite{I}.  \qed
\eex

\br
By abuse of notation, we will set, for any variety $X$ with a Cartier divisor $D$,
$$\chi_y(X,\cO_X(D)):=\chi_y(X,D;\bQ_X^H)=\chi(X;mC_y(X) \otimes \cO_X(D))$$
and, if $X$ is pure dimensional of complex dimension $n$, 
$$I\chi_y(X,\cO_X(D)):=\chi_y(X,D;{IC'}_X^H),$$
where ${IC'}_X^H:=IC_X^H[-n]$.
\er

We next prove a generalized Hirzebruch-Riemann-Roch for $(X,D)$ as above, with any coefficients $\cM \in D^b\MHM(X)$. 

\bt\label{gHRR}(generalized Hirzebruch-Riemann-Roch with arbitrary coefficients)\newline
Let $X$ be a compact complex  algebraic variety, let $D$ be a  Cartier divisor on $X$, and fix $\cM \in D^b\MHM(X)$. Then the Hodge polynomial of $(X,D,\cM)$ is computed by the formula:
\be\label{Hp}
\chi_y(X,D;\cM)=\int_{X} ch(\cO_X(D)) \cap T_{y*}([\cM]).
\ee
\et

\begin{proof}
This is an application of the module property of the Todd class transformation, according to which, 
if $\alpha \in K_0(X)$ and $\beta \in K^0(X)$, then 
\be\label{m}
td_*(\beta \otimes \alpha)=ch(\beta) \cap td_*(\alpha),
\ee
with $ch$ denoting the {\it Chern character}.
We therefore have the following sequence of equalities:
\be
\begin{split}
\int_{X} ch(\cO_X(D)) \cap T_{y*}([\cM])&\overset{}{=}\sum_{p\in \bZ} \big[ \int_{X} ch(\cO_X(D)) \cap td_*([Gr^F_{-p}\DR(\cM)]) \big] \cdot (-y)^p\\
&\overset{(\ref{m})}{=}\sum_{p\in \bZ} \big[ \int_{X}  td_*([Gr^F_{-p}\DR(\cM) \otimes \cO_X(D) ]) \big] \cdot (-y)^p\\
&\overset{(\ast)}{=}\sum _{p\in \bZ} \chi\left(X,Gr^F_{-p}\DR(\cM) \otimes \cO_X(D) \right) \cdot (-y)^p\\
&= \chi(X; \DR_y([\cM]) \otimes \cO_X(D)) \\
&=\chi_y(X,D;\cM),
\end{split}
\ee
where ($\ast$) follows from the degree property of the Todd class transformation.
\end{proof}

\bc
If $D=0$, we get 
\be\label{zero}
\chi_y(X,0;\cM)=\int_{X}  T_{y*}([\cM]) =\chi_y(X;\cM),
\ee
which is the usual Hodge polynomial of $\cM$ as in \eqref{eq2}. 
\ec

We next prove a twisted version of the duality formula \eqref{eq5}.
\bc[Duality formula]
Let $X$ be a compact complex  algebraic variety, let $D$ be a  Cartier divisor on $X$, and fix $\cM \in D^b\MHM(X)$. Then 
\be\label{duf}
\chi_y(X,D;\bD_X\cM)=\chi_{1/y}(X,-D;\cM).
\ee
\ec

\begin{proof}
Since Grothendieck duality on a point space is the identity, the assertion follows by pushing down to a point the following duality formula for twisted Hodge-Chern classes:
\begin{equation*}
\begin{split}
\DR_y([\bD_X\cM]) \otimes \cO_X(D)&= \bD \left( \DR_{1/y}([\cM]) \right) \otimes \cO_X(D)\\ &=\bD \left( \DR_{1/y}([\cM]) \otimes \cO_X(-D) \right).
\end{split}
\end{equation*}
Here, the first equality follows from Remark \ref{r9}, and the second by the module property of Grothendieck duality, i.e., if $\cL$ a line bundle and $\cF$ is a coherent sheaf of $\cO_X$-modules on $X$, then $\bD(\cL \otimes \cF)=\cL^\vee \otimes \bD(\cF)$, with $\cL^\vee$ the dual line bundle; e.g., see \cite[Part 1, Sect.7]{FM}. 
\end{proof}

\bex
If $\cM$ is a self-dual pure Hodge module of weight $k$, hence $\bD_X(\cM)\simeq \cM(k)$, then
\be\label{dual2} \chi_y(X,D;\cM)=(-y)^k \cdot \chi_{1/y}(X,-D;\cM).\ee
This also applies to shifted pure Hodge modules since for any mixed Hodge module complex, $\bD_X(\cM[m])\simeq (\bD_X\cM)[-m]$, and 
even shifts do not change the class in the Grothendieck group.
For example, $IC_X^H$ is a pure Hodge module of weight $n$, for $X$ pure $n$-dimensional, so formula \eqref{dual2} applies to ${IC'}_X^H$ to give:  $I\chi_y(X,\cO_X(D))=(-y)^n \cdot I\chi_{1/y}(X,\cO_X(-D))$.\qed
\eex

More explicit formulae for $\chi_y(X,D;\cM)$ can be obtained via the projection formula by expressing $T_{y*}([\cM])$ in terms of a stratification of $X$, as in Propositions \ref{p7} and \ref{p8}. As an application, we'll work this out in detail in the next section in the toric context.


\section{Characteristic class formulae in toric geometry}\label{sec5}

In this section we apply the above formulae in the case of toric varieties (e.g., associated to a full-dimensional lattice polytope), which are known to be Whitney stratified by the torus orbits. We are also interested here in the classes $IT_{y*}$ of toric varieties, and what they ``count'' when the variety comes from a polytope.

In view of formula \eqref{f5}, the key additional input in the toric situation is the following result derived in \cite{MS}:
\bp\cite[Prop.3.2]{MS}
Let $X_{\Sigma}$ be the toric variety defined by the fan $\Sigma$. For any cone $\sigma \in \Sigma$, with orbit $O_{\sigma}$ and inclusion $i_{\sigma} : O_{\sigma} \hookrightarrow \overline{O}_{\sigma}=V_{\sigma}$ in the orbit closure, we have:
\be\label{drt}
mC_y([i_\sig])=\DR_y(\big[(i_{\sigma})_!\bQ^H_{O_{\sigma}}\big])=(1+y)^{\dim(O_{\sigma})} \cdot [\omega_{V_{\sigma}}],
\ee
where $\omega_{V_{\sigma}}$ is the canonical sheaf on the toric variety $V_{\sigma}$.
\ep

This follows from Proposition \ref{p8} by using the fact that $\Omega_{Z}^1(\log(D))$ is the trivial sheaf of rank $\dim(O_\sig)$, together with $R(\pi_Z)_*\cO_Z(-D)\simeq \omega_{V_\sig}$ using a toric resolution $\pi_Z:Z \to V_\sig$ in Proposition \ref{p8}.

\br
In the above notations, we get from Remark \ref{ps} the following formula, using the fact that $V_{\sigma}$ has rational singularities, so $R(\pi_Z)_*\cO_Z\simeq \cO_{V_\sig}$, 
\be\label{drtc}
\DR_y(\big[(i_{\sigma})_*\bQ^H_{O_{\sigma}}\big])=(1+y)^{\dim(O_{\sigma})} \cdot [\cO_{V_{\sigma}}].
\ee
\er

Together with formula \eqref{f5}, this gives the following:
\bt\label{t28}
Let $X_{\Sigma}$ be the toric variety defined by the fan $\Sigma$. For each cone $\sigma \in \Sigma$ with orbit $O_{\sigma}$, let $k_{\sigma} : V_{\sigma} \hookrightarrow X_{\Sigma}$ be the inclusion of the orbit closure, and fix a point $x_\sigma \in O_{\sigma}$ with inclusion $i_{x_\sig}: \{x_\sig\} \hookrightarrow X_\Sig$.  Let $\cM \in D^b\MHM(X_\Sigma)$ (resp., $[\cM]\in K_0(\MHM(X))$)  be a (Grothendieck class of a) mixed Hodge module complex on $X_\Sigma$ which is (virtually) stratum-wise constant with respect to the stratification given by the torus orbits $O_{\sigma}$, $\sig \in \Sig$. 
Then:
\be\label{drt2g}
\DR_y([\cM])=\sum_{\sigma \in \Sigma} \, \chi_y(i_{x_\sig}^*[\cM]) \cdot (1+y)^{\dim(O_{\sigma})} \cdot (k_\sigma)_* [\omega_{V_{\sigma}}].
\ee
\be\label{tytg}
\begin{split}
T_{y*}([\cM])&=\sum_{\sigma \in \Sigma}\, \chi_y(i_{x_\sig}^*[\cM]) \cdot (1+y)^{\dim(O_{\sigma})} \cdot (k_\sigma)_* td_*([\omega_{V_{\sigma}}])\\
&=\sum_{\sig \in \Sig}\, \chi_y(i_{x_\sig}^*[\cM]) \cdot (-1-y)^{\dim(O_{\sig})} \cdot 
(k_{\sig})_*\left( td_*(V_{\sig})\right)^{\vee}.
\end{split}
\ee
In particular, if $X$ is compact, 
$$
\chi_y(X;\cM) =\sum _{\sig \in \Sig}\, \chi_y(i_{x_\sig}^*[\cM]) \cdot  (-1-y)^{\dim(O_{\sig})}.
$$

\et 

For the second equality in \eqref{tytg}, we use the duality property of $td_*$, namely, if $X$ is a Cohen-Macaulay variety (e.g., a toric variety like the closure of a torus orbit, see \cite[Thm.9.2.9]{CLS}), then we have (e.g., see \cite[Ex.18.3.19]{Fu3}):
\be\label{5}
\bD([\cO_X])=(-1)^{\dim(X)} \cdot [\omega_X], \ \ \ td_k(X)=(-1)^{\dim(X)-k} td_k ([\omega_X]).
\ee

\bex[Toric fibrations]\label{extf}
Let $f:Y \to X$ be a proper toric morphism of toric varieties, with the corresponding lattice homomorphism $f_N:N' \to N$ surjective (i.e., $f$ is a toric fibration  in the sense of \cite[Prop.2.1]{CMM}). 
Let $\Sig'$, $\Sig$ be the fans of $Y$, resp., $X$. Since $f$ is a toric fibration, a $\bT'$-orbit $O_{\sig'}$ ($\sig' \in \Sig'$) is mapped by $f$ to a $\bT$-orbit $f(O_{\sig'})=O_{\sig}$ ($\sig \in \Sig$), such that the restriction map $f_{\sig'}=f\vert_{O_{\sig'}}:O_{\sig'} \to O_{\sig}$ is isomorphic to a projection $O_{\sig'} \simeq O_{\sig} \times (\bC^*)^\ell \to O_{\sig}$ and $\ell=\dim(O_{\sig'}) - \dim(O_{\sig})$ the relative dimension of $f_{\sig'}$ (see \cite[Lem.2.6 and Prop.2.7]{CMM}). Let
\be\label{dl} d_\ell(Y/\sig):=\vert \Sigma_\ell(Y/\sig) \vert \ee with
\be\label{dl2} \Sigma_\ell(Y/\sig):=\{\sig' \in \Sig' \mid  f(O_{\sig'})=O_{\sig}, \ \ell=\dim(O_{\sig'}) - \dim(O_{\sig})\}.\ee
Under the above notations and assumptions, we have
\be\label{fn90}
f_* mC_y(Y)=\sum_{\sig \in \Sig}\chi_y (f^{-1}(s_\sig))  \cdot (1+y)^{\dim(O_\sig)} \cdot (k_\sig)_*[\omega_{V_{\sig}}] 
\ee
for $s_\sig \in O_\sig$, with 
$$\chi_y (f^{-1}(s_\sig))=\sum_{\ell \geq 0} d_\ell(Y/\sig) \cdot (-y-1)^\ell.$$
A similar formula holds for $f_* T_{y*}(Y)$. See \cite[Prop.3.9]{CMSSEM} for an equivariant analogue of these formulae. \qed
\eex

Let us next recall that an open affine subvariety $U_\sig \subset X_\Sig$ corresponding to a cone $\sig \in \Sig$, is of global product type $$U_\sig=Z_\sig \times O_\sig,$$ as follows. Choose a splitting $N=N_\sig \oplus {N'}$ of $N$, where $N_\sig=N \cap \langle \sig \rangle$ is the lattice spanned by $\sig$, with a corresponding splitting $\bT=\bT_\sig \times \bT'$ of $\bT$. Here $Z_\sig$ is the $\bT_\sig$-toric variety of $\sig$ in $N_\sig \otimes \bR$, and $\bT'\simeq O_\sig$.

\bl\label{vcc} Let $X=X_\Sig$ be the toric variety with torus $\bT$ associated to a fan $\Sigma$ in $N_\bR=N\otimes \bR$, with the Whitney stratification $\cS$ given by the $\bT$-orbits. Then:
\begin{itemize}
\item[(a)] $({j_S})_*\bQ_S^H$ (resp., $[({j_S})_*\bQ_S^H]\in K_0(\MHM(X))$)  is (virtually) stratum-wise (Tate-type) constant for each $S \in \cS$, with $j_S:S \hookrightarrow X$ the inclusion. 
\item[(b)] $IC_X^H$ ((resp., $[IC_X^H] \in K_0(\MHM(X))$)  is (virtually) stratum-wise (Tate-type) constant along $\cS$.
\item[(c)]
If $\cM$ (resp., $[\cM]$) is (virtually) stratum-wise (Tate-type) constant along $\cS$, 
then so is its dual $\bD_X \cM$ (resp., $\bD_X [\cM]$).
\end{itemize}
\el

\begin{proof}
It suffices to prove the above properties on each open affine subvariety $U_\sig \subset X$, for $\sig \in \Sig$, with the splitting $U_\sig=Z_\sig \times O_\sig,$ as recalled above. Let $p_\sig:U_\sig \to Z_\sig$ be the smooth projection map with fiber $O_\sig$. Then $(a)$ follows from the smooth base change property for the underlying constructible complex. Similarly, smooth pullback of mixed Hodge module complexes preserves the $K_0(\MHS^p)^{\rm Tate}$-module structure. Item $(b)$ can also be checked on the underlying perverse sheaf $IC_X$, where it follows by induction using $(a)$ and the Deligne construction. Again, these operations preserve the $K_0(\MHS^p)^{\rm Tate}$-module structure (see also \cite[(1.7.6)]{Sa}). For $(c)$, the stratum-wise constant case follows inductively from the fact that duality commutes (up to a shift) with the smooth pullback. Applying this to $[(j_S)_!\bQ^H_X]$ as in Example \ref{exvcc},  one proves the virtual stratum-wise constant case since duality commutes with cross-products (and duality on a point preserves the Tate type).
\end{proof}

\br In the notations of the Introduction, 
the subgroup $$\chi_{Hdg}\left(\tau(F^{\bT}(X)[y^{\pm 1}])\right) \subset K_0(\MHM(X))$$ corresponds to the subgroup given by classes $[\cM]$ wich are virtually stratum-wise Tate-type constant along the $\bT$-orbit stratification of the toric variety $X$.
\er

Theorem \ref{t28} and Lemma \ref{vcc}(b) yield the following (the case of the constant Hodge module was already considered by the authors in \cite{MS}): 
\bc\label{c28}
Let $X_{\Sigma}$ be the toric variety defined by the fan $\Sigma$. For each cone $\sigma \in \Sigma$ with orbit $O_{\sigma}$, let $k_{\sigma} : V_{\sigma} \hookrightarrow X_{\Sigma}$ be the inclusion of the orbit closure, and fix a point $x_\sigma \in O_{\sigma}$.   
Then we have the following formulae:
\be\label{drt2}
\DR_y([{IC'}_{X_{\Sigma}}^H])=\sum_{\sigma \in \Sigma} \, \chi_y({IC'}_{X_{\Sigma}}^H\vert_{x_\sigma}) \cdot (1+y)^{\dim(O_{\sigma})} \cdot (k_\sigma)_* [\omega_{V_{\sigma}}].
\ee
\be\label{tyt}
\begin{split}
IT_{y*}(X_{\Sigma})&=\sum_{\sigma \in \Sigma}\, \chi_y({IC'}_{X_{\Sigma}}^H\vert_{x_\sigma}) \cdot (1+y)^{\dim(O_{\sigma})} \cdot (k_\sigma)_* td_*([\omega_{V_{\sigma}}])\\
&=\sum_{\sigma \in \Sigma}\, \chi_y({IC'}_{X_{\Sigma}}^H\vert_{x_\sigma}) \cdot  (-1-y)^{\dim(O_{\sig})} \cdot 
(k_{\sig})_*\left( td_*(V_{\sig})\right)^{\vee}. 
\end{split}
\ee
In particular, if $X_\Sig$ is a simplicial toric variety then $\chi_y({IC'}_{X_{\Sigma}}^H\vert_{x_\sigma})=\chi_y(\bQ_{X_{\Sigma}}^H\vert_{x_\sigma})=1$ for all $\sig \in \Sig$.
\ec

\bex[$y=1$, signature]
The intersection homology signature of a projective toric variety $X_\Sigma$ is computed by
\be\label{sig}
sign(X_\Sigma)=\sum_{\sigma \in \Sigma}\, \chi_1({IC'}_{X_{\Sigma}}^H\vert_{x_\sigma}) \cdot  (-2)^{\dim(O_{\sig})}.
\ee
For the case of a simplicial toric variety, see e.g., \cite[formula (1.12)]{MS}.\qed
\eex


\section{Applications to weighted lattice point counting}\label{sec6}

Let us recall how the classical lattice point counting in a polytope relates to characteristic classes.
Let $M\cong \bZ^n$ be a lattice and $P \subset M_{\bR} \cong \bR^n$ be a full-dimensional lattice polytope with associated projective toric variety $X_P$, inner normal fan $\Sig_P$ and ample Cartier divisor $D_P$.
By the classical work of Danilov \cite{Dan},
the (dual) Todd classes of $X_P$ can be used for counting the number of lattice points in (the interior of) a lattice polytope $P$, namely,
\be\label{D1}
| P \cap  M | =\int_{X_P} ch(\cO_{X_P}(D_P)) \cap td_*(X_P) =\sum_{k\geq 0} \frac{1}{k!} \int_{X_P} [D_P]^k \cap td_k(X_P), 
\ee
\be\label{D2}
| \Int({P}) \cap M | =\int_{X_P} ch(\cO_{X_P}(D_P)) \cap td_*([\omega_{X_P}])=\sum_{k\geq 0} \frac{(-1)^n}{k!} \int_{X_P} [-D_P]^k \cap td_k(X_P),
\ee
where the summation  is over the faces $Q$ of $P$, $|-|$ denotes the cardinality, and $\Int(P)$ is the interior of $P$.
This also gives the coefficients of the {\it Ehrhart polynomial of $P$} counting the number of lattice points in the dilated polytope $\ell P:=\{\ell \cdot u \ | \ u \in P \}$ for a positive integer $\ell$:
$${\rm Ehr}_P(\ell):=|\ell P \cap M|=\sum_{k=0}^d a_k \ell^k,$$ with $$a_k= \frac{1}{k!}  \int_{X_P}  [D_P]^k \cap td_k(X_P).$$
Note that by the above formulae one gets ${\rm Ehr}_P(0)=a_0=\int_{X_P} td_*(X_P) =1$, and the reciprocity identity 
\be\label{cld} {\rm Ehr}_P(-\ell)=(-1)^{n}\cdot | \Int({\ell P}) \cap M |.\ee

Recall that in \cite{MS} the authors showed that the  homology Hirzebruch class 
\be\label{tytc}
\begin{split}
T_{y*}(X_P)&=\sum_{\sigma \in \Sigma_P}\,  (1+y)^{\dim(O_{\sigma})} \cdot (k_\sigma)_* td_*([\omega_{V_{\sigma}}])\\
\end{split}
\ee
can be used to count lattice points in a lattice polytope $P$ so that points in the interior of a face $Q$ of $P$ carry the weight $(1+y)^{\dim(Q)}$. More precisely, one has (see \cite[Theorem 1.3]{MS}):
\be\label{i5}
\sum_{Q \preceq P} (1+y)^{\dim(Q)} \cdot | \Relint(Q) \cap M |=\int_{X} ch(\cO_{X_P}(D_P)) \cap T_{y*}(X_P)=\chi_y(X_P, \cO_{X_P}(D_P)),
\ee
with $\Relint(Q)$ the relative interior of a face $Q$.

The classes $T_{y*}([\cM])$, and in particular $IT_{y*}(X_{P})$, serve a similar purpose, but the weights need to be adjusted, as the formulae of Theorem \ref{t28} suggest. 

\bp[Weighted lattice point counting]\label{wlpca}
If $X=X_P$ is the projective toric variety with ample Cartier divisor $D=D_P$ associated to a full-dimensional lattice polytope $P\subset M_{\bR} \cong \bR^n$, and $\ell \in \bZ_{>0}$, then for any 
$\cM \in D^b\MHM(X)$ (resp., $[\cM]\in K_0(\MHM(X))$) a (Grothendieck class of a) mixed Hodge module complex on $X$ which is (virtually) stratum-wise constant with respect to the stratification given by the torus orbits, we have 

\be\label{wca}
\begin{split}
\chi_y(X,\ell D;\cM) &= \int_{X} e^{\ell[D]} \cap T_{y*}([\cM])\\
&=
\sum_{k=0}^{n} \left(\frac{1}{k!} \int_{X} [D]^k \cap T_{y*}([\cM])  \right) \ell^k \\
 &=\sum_{Q \preceq P} \chi_y(i_{x_Q}^*[\cM]) \cdot (1+y)^{\dim(Q)} \cdot | \Relint({{\ell}} Q) \cap M | , 
\end{split}
\ee
where the summation on the right is over the faces $Q$ of $P$, and $x_Q\in O_{\sigma_Q}$ is a point in the orbit of $X_P$ associated to (the cone $\sigma_Q$ in the inner normal fan of $P$ corresponding to) the face $Q$.
\ep

\begin{proof} For a face $Q$ of $P$, denote by $k_Q:V_{\sig_Q}:=X_Q \hookrightarrow X$ the inclusion of the orbit closure associated to the (cone $\sigma_Q$ of the) face $Q$. Note that we have $\dim_\bR(Q)=\dim_\bR(O_{\sig_Q}).$ Then by \eqref{tytg} the following  equality holds:
$$\int_{X} ch(\cO_{X}(\ell D)) \cap T_{y*}([\cM])
= \sum_{Q \preceq P} \chi_y(i_{x_Q}^*[\cM]) \cdot  (1+y)^{\dim(Q)} \int_{X} ch(\cO_{X}(\ell D)) \cap (k_Q)_*td_*([\omega_{X_Q}]).$$

It remains to prove that for any face $Q$ of $P$, we have that:
\be\label{lem2}
 \int_{X} ch(\cO_{X}(\ell D)) \cap (k_Q)_*td_*([\omega_{X_Q}]) = | \Relint(\ell Q) \cap M |.
\ee
This follows exactly as in \cite[Theorem 1.3]{MS} by using Danilov's formula \eqref{D2} for the orbit closure $X_Q$, together with  
\be\label{lem2}
 \int_{X} ch(\cO_{X}(\ell D)) \cap (k_Q)_*td_*([\omega_{X_Q}]) 
 = \int_{X_Q}  ch(\cO_{X_Q}(\ell D_Q)) \cap td_*([\omega_{X_Q}]). 
\ee
The latter equality uses the projection formula, together with $(k_Q)^*(\cO_{X_P}(\ell D_P)=\cO_{X_Q}(\ell D_Q)$.
\end{proof}

\bd[Weighted Ehrhart polynomials]
In the above notations, the expression
$$E_{P,\cM}(\ell, y):=\sum_{Q \preceq P} \chi_y(i_{x_Q}^*[\cM]) \cdot (1+y)^{\dim(Q)} \cdot | \Relint({{\ell}} Q) \cap M | \in \bZ[y^{\pm 1}]$$
is obtained by evaluating the polynomial 
$$E_{P,\cM}(z, y)=\sum_{k=0}^{n} \left(\frac{1}{k!} \int_{X_P} [D_P]^k \cap T_{y*}([\cM])  \right) z^k$$
at $z=\ell$. We call this polynomial the {\it weighted Ehrhart polynomial} of $P$ with coefficients induced from 
$\cM \in D^b\MHM(X_P)$ (resp., $[\cM]\in K_0(\MHM(X_P))$), a (Grothendieck class of a) mixed Hodge module complex on $X_P$ which is (virtually) stratum-wise constant with respect to the stratification given by the torus orbits.

More generally, for a Laurent polynomial {\it weight function} $$f:\{\text{faces of } P\} \to \bZ[y^{\pm 1}], \ \ \ Q \mapsto f_Q(y)$$ defined on the set of faces of $P$, we define as in the Introduction an associated weighted Ehrhart polynomial by
$$E_{P,f}(\ell, y):=\sum_{Q \preceq P} f_Q(y) \cdot (1+y)^{\dim(Q)} \cdot | \Relint({{\ell}} Q) \cap M | .$$
\ed

\br\label{r34}
Note that any $\cM \in D^b\MHM(X_P)$ (resp., $[\cM]\in K_0(\MHM(X_P))$), a (Grothendieck class of a) mixed Hodge module complex on $X_P$ which is (virtually) stratum-wise constant with respect to the stratification given by the torus orbits, 
induces a weight function on the faces of $P$ defined by $f_Q(y)=\chi_y(i_{x_Q}^*[\cM])$, with $x_Q\in O_{\sigma_Q}$ a point in the orbit of $X_P$ associated to (the cone $\sigma_Q$ in the inner normal fan of $P$ corresponding to) the face $Q$ and $i_{x_Q}: \{x_Q\} \to X_P$ the inclusion map. Moreover, any weight function $f$ can be obtained in this way from such a mixed Hodge module complex $\cM$, e.g., by using direct sums and shifts of $(j_Q)_! \bQ^H_{O_{\sig_Q}}(n_Q)$, with $n_Q \in \bZ$ and $j_Q:O_{\sig_Q}\hookrightarrow X$ the inclusion (or, using the subgroup $\chi_{Hdg}\left(\tau(F^{\bT}(X)[y^{\pm 1}])\right)$ as in the Introduction). In particular, the expression $E_{P,f}(\ell, y)$ is obtained by evaluating a {\it polynomial} $E_{P,f}(z, y)$ at $z=\ell$.
\er

By using the first equality of \eqref{wca}, we obtain the following:
\bp
In the above notations, 
\be
E_{P,\cM}(0, y)=\chi_y(X_P,0;\cM)=\chi_y(X_P;\cM)
\ee
is the usual Hodge polynomial of $\cM$.
\ep

\bex\label{ex65}
Let us next give some examples for Proposition \ref{wlpca}, formulated in terms of Ehrhart theory. If $X=X_P$ is the projective toric variety with inner normal fan $\Sig=\Sig_P$ and ample Cartier divisor $D=D_P$ associated to a full-dimensional lattice polytope $P\subset M_{\bR} \cong \bR^n$, and $\ell \in \bZ_{>0}$, for a face $Q$ of $P$, denote by $j_Q:O_{\sig_Q} \hookrightarrow X$ the inclusion of the orbit associated to the (cone of the) face $Q$.
\begin{enumerate}
\item[(a)] For $\cM=(j_Q)_! \bQ^H_{O_{\sig_Q}}$, we get
\be\label{69} E_{P,\cM}(\ell, y)=(1+y)^{\dim(Q)} \cdot | \Relint({{\ell}} Q) \cap M | \in \bZ[y],\ee
with $E_{P,\cM}(0, y)=\chi_y(O_{\sig_Q})=(-y-1)^{\dim(Q)}$.

In particular, by evaluating $E_{P,\cM}(\ell, y)$ at $y=0$, we get the polynomial
$E_{P,\cM}(\ell, 0)= |\Relint({{\ell}} Q) \cap M |$ counting lattice points in the relative interior of the dilated face $\ell Q$, with 
$E_{P,\cM}(0, 0)=(-1)^{\dim(Q)}$.

\item[(b)] Let $X':=X_{P'}$ be a torus-invariant closed algebraic subset of $X=X_P$ corresponding to a polytopal subcomplex $P' \subseteq P$ (i.e., a closed union of faces of $P$), with $\cM=\bQ_{X'}^H$ pushed forward to the ambient $X_P$. Then 
\be
E_{P,\cM}(\ell, y)= \sum_{Q \preceq P'} (1+y)^{\dim(Q)} \cdot | \Relint(\ell Q) \cap M |  \in \bZ[y],
\ee
where the summation is over the faces $Q$ of $P'$. Hence, 
 $E_{P,\cM}(\ell, 0)=|\ell P' \cap M|$, so if $y=0$ we recover the {\it classical Ehrhart polynomial} ${\rm Ehr}_{P'}(\ell)$ for $P'$, and 
$$E_{P,\cM}(0, y)=\chi_y(X')=\sum_{Q \preceq P'} (-y-1)^{\dim(Q)},$$
which gives that $$E_{P,\cM}(0, 0)=\sum_{Q \preceq P'} (-1)^{\dim(Q)}=\chi(P')$$ is the topological Euler characteristic of $P'$; see also \cite[formula (4.4)]{MS}. Similarly, for the inclusion of the open complement $j_U:U=X \setminus X' \hookrightarrow X$, with $\cM=(j_U)_!\bQ^H_U$, we get similar formulae summing over faces $Q$ of $P$ not contained in $P'$.

\item[(c)] For $\cM=(j_Q)_* \bQ^H_{O_{\sig_Q}}$, we get by \eqref{drtc} and Danilov's formula \eqref{D1} that 
\be\label{pu} E_{P,\cM}(\ell, y)=(1+y)^{\dim(Q)} \cdot | {{\ell}} Q \cap M | \in \bZ[y],\ee
with $E_{P,\cM}(0, y)=\chi_y(V_{\sig_Q})$, $E_{P,\cM}(\ell, 0)=| {{\ell}} Q \cap M |$ fitting also for $y=0$ with the classical Ehrhart polynomial for $Q$, and $E_{P,\cM}(0, 0)=\chi(Q)$. Note that our choice of $\cM$ satisfies the stratum-wise constant  assumption, e.g., by Lemma \ref{vcc}.
\qed
\end{enumerate}
\eex


Let $X=X_P$ be as before the projective toric variety with inner normal fan $\Sig=\Sig_P$ and ample Cartier divisor $D=D_P$ associated to a full-dimensional lattice polytope $P\subset M_{\bR} \cong \bR^n$.
Then another important example is provided by $\cM={IC'}^H_{X_P}$.
As explained in the Introduction, in the case when $P$ contains the origin in its interior (so that its polar polytope $P^\circ$ is defined), the stalk contributions $\chi_y({IC'}_{X_P}^H\vert_{x_Q})$ have the following combinatorial description (see, e.g., \cite{BM,F, DL,Sa}, etc):
\be\label{co}
\chi_y({IC'}_{X_P}^H\vert_{x_Q})=\wti{g}_Q(-y):=g_{Q^\circ}(-y),
\ee
with $Q$ a face of $P$ and $Q^\circ$  the corresponding face of the polar polytope $P^\circ$, and $g_{Q^\circ}(t) \in \bZ[t]$ the Stanley {\it $g$-polynomial} of $Q^\circ$, cf. \cite[Sect.3]{St}.
We then get:


\bc\label{wlpc}
If $X=X_P$ is the projective toric variety with ample Cartier divisor $D=D_P$ associated to a full-dimensional lattice polytope $P\subset M_{\bR} \cong \bR^n$, and $\ell \in \bZ_{>0}$, then:
\be\label{wc}
\begin{split}
\int_{X} e^{\ell[D]} \cap IT_{y*}(X)&=
I\chi_y(X,\cO_X(\ell D)) \\ &=\sum_{Q \preceq P} \chi_y({IC'}_{X}^H\vert_{x_Q}) \cdot (1+y)^{\dim(Q)} \cdot |\Relint({{\ell}} Q) \cap M |  \\
&=E_{P,IC'}(\ell, y) \in \bZ[y], 
\end{split}
\ee
where the summation is over the faces $Q$ of $P$, and $x_Q\in O_{\sigma_Q}$ is a point in the orbit of $X=X_P$ associated to (the cone $\sigma_Q$ in the inner normal fan of $P$ corresponding to) the face $Q$. In particular,
$E_{P,IC'}(0, y)=I\chi_y(X)$, with $E_{P,IC'}(0, 0)=I\chi_0(X)$. Moreover, for $y=1$, one gets a weighted Ehrhart polynomial
\be\label{s2}
E_{P,IC'}(\ell, 1)=\sum_{Q \preceq P} \chi_1({IC'}_{X_P}^H\vert_{x_Q}) \cdot 2^{\dim(Q)} \cdot |\Relint({{\ell}} Q) \cap M |
\ee
whose constant term $E_{P,IC'}(0, 1)=I\chi_1(X)=sign(X)$ is the intersection cohomology signature of the projective toric variety $X=X_P$. 
\ec

\bex
Assuming in addition that the origin is an interior point of $P$, the formulas of the above corollary have a combinatorial interpretation, since the weight function for $\cM={IC'}^H_{X_P}$ is $f_Q(y)=\wti{g}_Q(-y)$. With this weight function, we get the corresponding weighted Ehrhart polynomial 
\be E_{P,f}(\ell, y)=\sum_{Q \preceq P} \wti{g}_Q(-y) \cdot (1+y)^{\dim(Q)} \cdot |\Relint({{\ell}} Q) \cap M | \in \bZ[y],
\ee
with \be\label{icco} I\chi_y(X)=E_{P,f}(0, y)=\sum_{Q \preceq P} \wti{g}_Q(-y) \cdot (-1-y)^{\dim(Q)},\ee hence for $y=1$ one gets via \eqref{hin}  the following combinatorial formula for the intersection cohomology signature of $X=X_P$:
\be\label{sigc} sign(X)=E_{P,f}(0, 1)=\sum_{Q \preceq P} \wti{g}_Q(-1) \cdot (-2)^{\dim(Q)}.\ee \qed
\eex

\br[Simple polytopes]
 If the full-dimensional polytope $P$ is simple, then $X=X_P$ is a simplicial toric variety, hence ${IC'}_{X}^H \simeq \bQ^H_{X}$. So, in this case,  $\chi_y({IC'}_{X}^H\vert_{x_Q})=1$ and formula \eqref{s2} reduces to:
 $$E_{P,IC'}(\ell, 1)=\sum_{Q \preceq P} 2^{\dim(Q)} \cdot |\Relint({{\ell}} Q) \cap M |
,$$
with constant term $sign(X)=\sum_{Q \preceq P} (-2)^{\dim(Q)}$; see also \cite[formula (1.12)]{MS}.
 \er
 
 \br[Intersection cohomology Poincar\'e polynomials] 
Making $\ell=0$ in formula \eqref{wc} yields via Example \ref{ex65}(a):
\be\label{wc2}
\begin{split}
I\chi_y(X)=\int_{X} IT_{y*}(X)&=\sum_{Q \preceq P} \chi_y({IC'}_{X}^H\vert_{x_Q}) \cdot (-1-y)^{\dim(Q)} , 
\end{split}
\ee
computing the intersection cohomology Hodge polynomials of $X=X_P$ in terms of local contributions and the combinatorics of the polytope. Moreover, 
using the fact that both $[IH^*(X)]$ and $[{IC'}_{X_P}^H\vert_{x_Q}]$ are of Tate type in $K_0(\MHS)$, upon substituting  $y=-t^2$ in \eqref{wc2}, one recovers a formula of Fieseler \cite[Theorem 1.1]{F} for the intersection cohomology Poincar\'e polynomial of $X$, see also  \cite{DL}, \cite{Sa}, etc. \er

We next discuss the effect of duality on weighted Ehrhart polynomials, in the form of the following {\it reciprocity formula}. This takes the form of a {\it purity} statement for pure Hodge modules.
\bt[Reciprocity and Purity]\label{repr} In the above notations, for any 
$\cM \in D^b\MHM(X)$ (resp., $[\cM]\in K_0(\MHM(X))$) a (Grothendieck class of a) mixed Hodge module complex on $X$ which is (virtually) stratum-wise constant with respect to the stratification given by the torus orbits, we have
\be\label{r1}
\begin{split}
E_{P,\cM}(-\ell, y)=E_{P, \bD_X\cM}(\ell, \frac{1}{y}).
\end{split}
\ee
In particular, if $\cM$ is such a self-dual pure Hodge module of weight $n$ on $X=X_P$,
then the following purity property holds:
\be\label{r1b}
E_{P,\cM}(-\ell, y)=(-y)^n \cdot E_{P, \cM}(\ell, \frac{1}{y}).
\ee
More generally, for any weight function $f$ on the faces of $P$, we have
\be\label{r2}
E_{P,f}(-\ell, y)=\sum_{Q \preceq P} f_Q({y}) \cdot (1+y)^{\dim(Q)} \cdot \left( (-1)^{\dim(Q)} \cdot |{{\ell}} Q \cap M | \right).
\ee
\et

\begin{proof}
Formula \eqref{r1} follows immediately by combining \eqref{wca} with the duality formula \eqref{duf}. Indeed, with $D=D_P$ the associated Cartier divisor of $P$, 
$$E_{P,\cM}(-\ell, y)=\chi_y(X,-\ell D;\cM)=\chi_{1/y}(X,\ell D;\bD_X\cM)=E_{P, \bD_X\cM}(\ell, \frac{1}{y}).$$ 
Formula \eqref{r1b} follows similarly by using \eqref{dual2}.

Let now $f$ be a weight function as above, and choose a representative $[\cM] \in K_0(\MHM(X_P))$ which is virtually stratum-wise constant along the $\bT$-orbits, such that 
$E_{P,f}=E_{P,\cM}$. Then $\bD_X[\cM]$ is also virtually stratum-wise constant along the $\bT$-orbits by Lemma \ref{vcc}. 
By \eqref{drt2g}, we have
\be
\begin{split}
\DR_y([\cM])&=\sum_{Q \preceq P} \, \chi_y(i^*_{x_Q}[\cM]) \cdot (1+y)^{\dim(Q)} \cdot (k_Q)_* [\omega_{V_{\sigma_Q}}] \\
&=\sum_{Q \preceq P} \, f_Q(y) \cdot  \DR_y \left([(j_Q)_! \bQ^H_{O_{\sig_Q}}] \right),
\end{split}
\ee
with $f_Q(y)=\chi_y(i^*_{x_Q}[\cM])$. Moreover, since on the smooth stratum $S$,  $\bQ_S^H \simeq {IC'}^H_S$ is a pure Hodge module of weight $\dim(S)$, we get as in Remark \ref{r9} that
$$\DR_y \left(\bD_X \left([(j_Q)_! \bQ^H_{O_{\sig_Q}}]\right) \right)
= (-y)^{-\dim(Q)} \cdot   \DR_y \left([(j_Q)_* \bQ^H_{O_{\sig_Q}}]\right). $$
Since $\DR_y$ commutes with dualities, up to exchanging $y$ by $1/y$, we then have
$$\DR_{1/y}(\bD_X[\cM])=\sum_{Q \preceq P} \, f_Q({y}) \cdot  (-y)^{\dim(Q)} \cdot \DR_{1/y}\left([(j_Q)_* \bQ^H_{O_{\sig_Q}}]\right).$$
Then, with $X=X_P$, $D=D_P$, we get: 
\begin{equation*}
\begin{split}
E_{P,f}(-\ell, y)&=\chi_{1/y}(X,\ell D;\bD_X\cM)=\int_{X} e^{\ell D} \cap td_* (\DR_{1/y}(\bD_X[\cM])) \\
&=\sum_{Q \preceq P} \, f_Q({y}) \cdot  (-y)^{\dim(Q)} \cdot \int_{X} e^{\ell D} \cap T_{1/y,*} \left([(j_Q)_* \bQ^H_{O_{\sig_Q}}]\right)\\
&\overset{\eqref{pu}}{=}\sum_{Q \preceq P} \, f_Q({y}) \cdot  (-y)^{\dim(Q)} \cdot \left(1+\frac{1}{y}\right)^{\dim(Q)} \cdot |\ell Q \cap M| \\
&=\sum_{Q \preceq P} \, f_Q({y})  \cdot \left(1+{y}\right)^{\dim(Q)} \cdot  \left(  (-1)^{\dim(Q)} \cdot |\ell Q \cap M| \right).
\end{split}
\end{equation*}
\end{proof}


\bex
The reciprocity formula \eqref{r1} can be applied to all cases of Example \ref{ex65}. In the (dual) cases $(a)$ and $(c)$, for $\cM=(j_Q)_! \bQ^H_{O_{\sig_Q}}$ with $\bD_X [\cM]=[(j_Q)_* \bQ^H_{O_{\sig_Q}}(\dim(Q))]\in K_0(\MHM(X))$, formula \eqref{r1} specializes to the {\it classical Ehrhart reciprocity} \eqref{cld}  for $Q$. Indeed, 
$$
E_{P,\cM}(\ell, y)=(1+y)^{\dim(Q)} \cdot | \Relint({{\ell}} Q) \cap M |
$$
and 
$$
E_{P,\bD_X\cM}(\ell, \frac{1}{y})=(1+y)^{\dim(Q)} \cdot \left(  (-1)^{\dim(Q)} \cdot |\ell Q \cap M| \right).
$$ \qed
\eex

\bex[$IC$-reciprocity]\label{icre}
Assume $P$ contains the origin in its interior. Then for 
$\cM={IC'}_X^H$, corresponding to the weight function $f_Q(y)=\wti{g}_Q(-y)$, the weighted Ehrhart polynomial becomes:
\be\label{ep1} E_{P,f}(\ell, y)=\sum_{Q \preceq P} \wti{g}_Q(-y) \cdot (1+y)^{\dim(Q)} \cdot |\Relint({{\ell}} Q) \cap M | \in \bZ[y].\ee
Moreover, by \eqref{r2} or classical Ehrhart reciprocity, we get
\be E_{P,f}(-\ell, \frac{1}{y})=\sum_{Q \preceq P} \wti{g}_Q(-\frac{1}{y}) \cdot (-y)^{-\dim(Q)} \cdot (1+y)^{\dim(Q)} \cdot |\ell Q \cap M| .\ee
Hence using \eqref{r1b} for ${IC'}^H_X$,  we get 
\be\label{ep2} E_{P,f}(\ell, y)=\sum_{Q \preceq P} \wti{g}_Q(-\frac{1}{y}) \cdot (-y)^{n-\dim(Q)} \cdot (1+y)^{\dim(Q)} \cdot |\ell Q \cap M|.\ee \qed
\eex

\br
Our definition \eqref{ep1} for $E_{P,f}(\ell, y)$ equals $G_{\varphi=1}(\ell,y)$ from \cite[formula (14)]{BGM} for the constant function $\varphi=1$, and formula \eqref{ep2} matches the definition of $G_{\varphi=1}(\ell,y)$ given in \cite[formula (9)]{BGM}. Moreover, our general reciprocity formula \eqref{r1b} recovers the reciprocity formula \cite[Theorems 1.3 and 2.6]{BGM} for the constant function $\varphi=1$.
\er

As another consequence of \eqref{r1b} for $\cM=IC^H_{X_P}$ and $y=1$, we obtain
the following generalization of \cite[(4.5), (4.6)]{MS} (where the simplicial case was considered, with all $\wti{g}_Q(-1)=1$). 
\bc
For a full-dimensional lattice polytope $P\subset M_\bR \simeq \bR^n$ with toric variety $X=X_P$, one has the following identity:
\be\label{wcL}
\sum_{Q \preceq P} \chi_1({IC'}_{X}^H\vert_{x_Q}) \cdot \left(\frac{1}{2}\right)^{{\rm codim}(Q)} \cdot |\Relint(\ell Q) \cap M |=\sum_{Q \preceq P}  \chi_1({IC'}_{X}^H\vert_{x_Q}) \cdot \left(-\frac{1}{2}\right)^{{\rm codim}(Q)} \cdot |\ell Q \cap M |,
\ee
where for each face $Q$ of $P$ we choose a point $x_Q$ in the $\bT$-orbit of $X$ corresponding (via the inner normal fan of $P$) to $Q$. If, moreover, $P$ contains the origin in its interior, then a similar formula is obtained with $\chi_1({IC'}_{X}^H\vert_{x_Q})=\wti{g}_Q(-1)$.
\ec

\begin{proof} Assume for simplicity that $P$ contains the origin in its interior. Then \eqref{wcL} follows by setting $y=1$ in formulae \eqref{ep1} and \eqref{ep2}, and then dividing the resulting expressions by $2^n$ with $n=\dim(P)=\dim(X)$. Here, we use the weights $f_Q(1)=\chi_1({IC'}_{X}^H\vert_{x_Q})=\wti{g}_Q(-1)$.  
\end{proof}




\end{document}